\documentclass[oneside]{amsart}

\usepackage{geometry}

\usepackage[english]{babel}
\usepackage{amsmath, amssymb, amsthm}
\usepackage{tikz}
\usepackage{tikz-cd}
\usepackage{adjustbox}
\usepackage{enumitem}
\usepackage{mathdots}
\usepackage{aligned-overset}
\usepackage{stmaryrd}
\usepackage[foot]{amsaddr}

\usepackage{subfiles}
\usepackage{hyperref}
\usepackage[noabbrev,capitalize]{cleveref}

\DeclareMathOperator{\Span}{Span}

\DeclareMathOperator{\coker}{coker}

\DeclareMathOperator{\Hom}{Hom}

\DeclareMathOperator{\aff}{aff}

\DeclareMathOperator{\op}{op}

\DeclareMathOperator{\modu}{Mod}
\DeclareMathOperator{\id}{id}

\DeclareMathOperator{\im}{im}
\DeclareMathOperator{\Spec}{Spec}

\DeclareMathOperator{\Rep}{Rep}

\DeclareMathOperator{\Gal}{Gal}

\DeclareMathOperator{\IHom}{\mathcal{H}om}
\DeclareMathOperator{\RIHom}{R\mathcal{H}\! \mathit{om}}

\DeclareMathOperator{\dgMod}{dgMod}

\DeclareMathOperator{\Irr}{Irr}

\DeclareMathOperator{\pt}{pt}
\DeclareMathOperator{\Perv}{Perv}

\DeclareMathOperator{\Spr}{Spr}
\DeclareMathOperator{\perf}{perf}

\DeclareMathOperator{\Sh}{Sh}
\DeclareMathOperator{\et}{\acute{e}t}

\DeclareMathOperator{\uHom}{\underline{\Hom}}

\DeclareMathOperator{\RuHom}{R\underline{\Hom}}

\DeclareMathOperator{\For}{For}

\DeclareMathOperator{\Rlim}{Rlim}
\DeclareMathOperator{\cont}{cont}
\DeclareMathOperator{\Ab}{Ab}
\DeclareMathOperator{\colim}{colim}
\DeclareMathOperator{\proet}{pro\acute{e}t}
\DeclareMathOperator{\Sbf}{\bf{S}}

\newtheorem{Theorem}{Theorem}[section]
\newtheorem{Proposition}[Theorem]{Proposition}
\newtheorem{Assumption}[Theorem]{Assumption}
\newtheorem*{Proposition*}{Proposition}
\newtheorem{Lemma}[Theorem]{Lemma}
\newtheorem{Corollary}[Theorem]{Corollary}
\newtheorem*{Theorem*}{Theorem}
\theoremstyle{definition}
\newtheorem{Definition}[Theorem]{Definition}
\newtheorem{Remark}[Theorem]{Remark}

\newtheorem{Example}[Theorem]{Example}
\newtheorem*{Conjecture*}{Conjecture}

\author{Jonas Antor}
\title{Formality in the Deligne-Langlands correspondence}
\email{\href{mailto:jonas.antor@maths.ox.ac.uk}{jonas.antor@maths.ox.ac.uk}}
\date{\today}

\begin{document}
\begin{abstract}
The Deligne-Langlands correspondence parametrizes irreducible representations of the affine Hecke algebra $\mathcal{H}^{\aff}$ by certain perverse sheaves. We show that this can be lifted to an equivalence of triangulated categories. More precisely, we construct for each central character $\chi$ of $\mathcal{H}^{\aff}$ an equivalence of triangulated categories between a perfect derived category of dg-modules $D_{\perf}(\mathcal{H}^{\aff}/(\ker(\chi)) - \dgMod)$ and the triangulated category generated by the corresponding perverse sheaves. The main step in this construction is a formality result that we prove for a wide range of `Springer sheaves'.
\end{abstract}

\maketitle

\section{Introduction}
\subsection*{Motivation and main result}
Affine Hecke algebras play an important role in the representation theory of $p$-adic groups. In fact, for each split $p$-adic group $\textbf{G}$ there is a (specialized) affine Hecke algebra $\mathcal{H}^{\aff}_q$ such that the category of $\mathcal{H}^{\aff}_q$-modules is equivalent to the category $\Rep_I(\textbf{G})$ of smooth $\textbf{G}$-representations that are generated by their Iwahori-fixed vectors \cite{borel1976admissible}. The category $\Rep_I(\textbf{G})$ is precisely the principal Bernstein block of $\textbf{G}$ (i.e. the block in the category of smooth $\textbf{G}$-representations that contains the trivial representation) \cite{bernstein1984centre,casselman1980unramified,bushnell1998smooth}. Let $G$ be the complex reductive group whose root datum is dual to that of $\textbf{G}$ and assume that $G$ has simply connected derived subgroup. In \cite{kazhdan1987proof} Kazhdan and Lusztig proved the Deligne-Langlands correspondence which parametrizes the irreducible representations of the affine Hecke algebra $\mathcal{H}^{\aff}$ by geometric data on the group $G$. We now recall this correspondence in the language of \cite{chriss2009representation}:
\medskip

Any irreducible representation of $\mathcal{H}^{\aff}$ admits a central character. These central characters are parametrized by semisimple conjugacy classes in $G \times \mathbb{G}_m$. For each semisimple pair $(s,q) \in G \times \mathbb{G}_m $, we denote the corresponding central character by $\chi_{(s,q)}: Z(\mathcal{H}^{\aff} ) \rightarrow \mathbb{C} $. Then the irreducible representations of $\mathcal{H}^{\aff}$ with central character $\chi_{(s,q)}$ are precisely the irreducible representations of the truncated affine Hecke algebra
\begin{equation}\label{eq: algebraic definition of truncated affine Hecke algebra}
\mathcal{H}^{\aff}_{(s,q)} := \mathcal{H}^{\aff}/(\ker(\chi_{(s,q)})).
\end{equation}
This algebra has a geometric incarnation in the world of constructible sheaves. Let $\mathcal{N}$ be the nilpotent cone of $G$ and $\mu: \tilde{\mathcal{N}} \rightarrow \mathcal{N}$ the Springer resolution. These varieties come with a canonical $G \times \mathbb{G}_m$-action where $\mathbb{G}_m$ acts by scaling. Passing to $(s,q)$-fixed points we obtain a morphism
\begin{equation*}
\mu^{(s,q)} : \tilde{\mathcal{N}}^{(s,q)} \rightarrow \mathcal{N}^{(s,q)}.
\end{equation*}
The corresponding $(s,q)$-Springer sheaf is defined as
\begin{equation}\label{eq: introduction sq springer sheaf}
\Sbf^{(s,q)}:= (\mu^{(s,q)})_* \mathcal{C}_{\tilde{\mathcal{N}}^{(s,q)}} \in D^b_c(\mathcal{N}^{(s,q)})
\end{equation}
where $\mathcal{C}_{\tilde{\mathcal{N}}^{(s,q)}}$ is the constant (perverse) sheaf on $\tilde{\mathcal{N}}^{(s,q)}$. The following isomorphism gives the affine Hecke algebra at $\chi_{(s,q)}$ a geometric interpretation \cite{kazhdan1987proof,chriss2009representation}:
\begin{equation}\label{eq: intro geometric incarnation of truncated affine hecke algebra}
\mathcal{H}^{\aff}_{(s,q)} \cong \Hom^*(\Sbf^{(s,q)}, \Sbf^{(s,q)}).
\end{equation}
Note that this isomorphism induces a grading on the algebra $\mathcal{H}^{\aff}_{(s,q)}$ which was not visible in its algebraic definition in \eqref{eq: algebraic definition of truncated affine Hecke algebra}. By the decomposition theorem \cite{beilinson2018faisceaux}, we can write
\begin{equation*}
\Sbf^{(s,q)} \cong \bigoplus_{j=1}^n X_j [k_j]
\end{equation*}
for certain $G(s)$-equivariant simple perverse sheaves $X_j \in   \Perv_{G(s)}(\mathcal{N}^{(s,q)})$ and integers $k_j \in \mathbb{Z}$. From this, one can deduce that there is a bijection
\begin{equation*}
\Irr(\mathcal{H}^{\aff}_{(s,q)}) \overset{1:1}{\longleftrightarrow} \{X_1, ..., X_n\}
\end{equation*}
or equivalently an embedding
\begin{equation*}
\Irr(\mathcal{H}^{\aff}_{(s,q)}) \hookrightarrow \Perv_{G(s)}(\mathcal{N}^{(s,q)}).
\end{equation*}
The simple objects in $ \Perv_{G(s)}(\mathcal{N}^{(s,q)})$ are parametrized by geometric data on $\mathcal{N}^{(s,q)}$, namely by the set of irreducible equivariant local systems on the $G(s)$-orbits of $\mathcal{N}^{(s,q)}$. Hence, we can interpret $\mathcal{N}^{(s,q)}$ as a variety of \textit{Langlands parameters} associated to the pair $(s,q)$. If $q \in \mathbb{G}_m$ is not a root of unity, the local systems that correspond to elements of $\Irr(\mathcal{H}^{\aff}_{(s,q)})$ can be characterized more explicitly as those that appear in the cohomology of a certain Springer fiber (c.f. \cite[Proposition 8.1.14]{chriss2009representation}). This is known as the \textit{Deligne-Langlands correspondence}.
\medskip

Following the categorical Langlands philosophy, one would like to lift the Deligne-Langlands correspondence to an equivalence of (triangulated) categories. Our main result establishes such an equivalence for each central character $\chi_{(s,q)}$. We denote by
\begin{equation*}
D_{\Spr}(\mathcal{N}^{(s,q)}) := \langle X_1, ..., X_n \rangle_{\Delta}
\end{equation*}
the full triangulated subcategory of $D^b_c(\mathcal{N}^{(s,q)})$ generated by the simple constituents of the $(s,q)$-Springer sheaf $\Sbf^{(s,q)}$.
\begin{Theorem*}(\cref{Theorem: derived Deligne-Langlands correspondence})
There is an equivalence of triangulated categories
\begin{equation*}
 D_{\perf}(\mathcal{H}^{\aff}_{(s,q)}- \dgMod) \cong D_{\Spr}(\mathcal{N}^{(s,q)})^{\op}
\end{equation*}
which identifies $\Sbf^{(s,q)}$ with the free dg-module $\mathcal{H}^{\aff}_{(s,q)}$. Here we consider $\mathcal{H}^{\aff}_{(s,q)}$ as a dg-algebra with vanishing differential and grading induced by the $\Hom^*$-grading in \eqref{eq: intro geometric incarnation of truncated affine hecke algebra}.
\end{Theorem*}
\subsection*{Formality of Springer sheaves}
The theorem above will be a consequence of a formality result that we prove for a wide range of `Springer sheaves': Let $G$ be a reductive group over an algebraically closed field $\overline{\mathbb{F}}$ (not necessarily of characteristic $0$), $B \subset G$ a Borel subgroup, $V$ a $G$-representation and $\{V^i \subset V \mid i \in I \}$ a finite collection of $B$-stable subspaces. Then for each $i \in I$ we consider the `morphism of Springer type'
\begin{equation*}
\mu_i : G \times^B V^i \rightarrow V, \quad (g,v) \mapsto gv
\end{equation*}
and the associated `Springer sheaves'
\begin{align*}
\Sbf^i&:= (\mu_i)_* \mathcal{C}_{G \times^B V^i}  \in D^b_c(V)\\
\Sbf &:= \bigoplus_{i \in I} \Sbf^i \in D^b_c(V).
\end{align*}
We will see in \cref{Corollary: the sqSpringer resolution is of Springer type} that the $(s,q)$-Springer sheaf $\textbf{S}^{(s,q)}$ from \eqref{eq: introduction sq springer sheaf} is a special case of this construction. We define the `Springer category' $D_{\Spr}(V)$ to be the full triangulated subcategory of $D^b_c(V)$ generated by the summands of $\Sbf$. Our goal is to establish an equivalence of triangulated categories
\begin{equation*}
D_{\Spr}(V)^{\op} \cong D_{\perf}(\Hom^*(\Sbf, \Sbf)-\dgMod).
\end{equation*}
The derived category $D_{\Spr}(V)^{\op}$ has a natural dg-enhancement (see the remarks below). Hence, by standard dg-techniques, which will be reviewed in \cref{Section: Derived categories and dg-algebras}, there is an equivalence of triangulated categories
\begin{equation}\label{eq: intro dg morita equivalence}
D_{\Spr}(V)^{\op}  \cong D_{\perf}(R_{\Sbf}-\dgMod)
\end{equation}
for a dg-algebra $R_{\Sbf}$ with
\begin{equation*}
H^*(R_{\Sbf}) = \Hom^*(\Sbf, \Sbf).
\end{equation*}
Thus, we need to show that we can replace the dg-algebra $R_{\Sbf}$ in the equivalence \eqref{eq: intro dg morita equivalence} with its cohomology, i.e. we need to prove that $R_{\Sbf}$ is \textit{formal}. We will do this by closely following the arguments in \cite[Appendix A,B]{polishchuk2019semiorthogonal}. First, we prove formality in the case where $\overline{\mathbb{F}}$ is the algebraic closure of a finite field $\mathbb{F} = \mathbb{F}_q$. In this case, $\Hom^*(\Sbf, \Sbf)$ can be equipped with a canonical Frobenius action. By analyzing the associated Steinberg variety, we will prove the following strong purity result for this action.
\begin{Theorem*}(\cref{Cor: Purity of the Ext algebra})
The canonical Frobenius action on $\Hom^i(\mathbf{S}, \mathbf{S})$ is given by multiplication with $q^{\frac{i}{2}}$.
\end{Theorem*}
We will then use a general "purity implies formality" result to deduce that $R_{\Sbf}$ is formal (the idea that purity implies formality goes back to \cite{deligne1975real,deligne1980conjecture}). This completes the proof in the positive characteristic case. Formality in characteristic $0$ can be deduced from the positive characteristic case using the "De $\mathbb{F}$ à $\mathbb{C}$" technique from \cite[§6]{beilinson2018faisceaux}. This last step will be discussed in more detail in \cref{Section: De F a C}.
\subsection*{Some remarks on dg-enhancements}
There are several technical problems one has to deal with to obtain a general "purity implies formality" result as we use it: One difficulty that arises is that the constructible derived category $D^b_c(X) = D^b_c(X, \overline{\mathbb{Q}}_{\ell})$ is not a genuine derived category, at least in its standard definition \cite{deligne1980conjecture,beilinson2018faisceaux}. In particular, $D^b_c(X, \overline{\mathbb{Q}}_{\ell})$ is not naturally dg-enhanced in this setting. It is shown in \cite{polishchuk2019semiorthogonal} that this problem can be solved by working with a dg-enhancement of $D^b_c(X,\overline{\mathbb{Q}}_{\ell})$ arising from an alternative constructions of the constructible derived category with integral coefficients $D^b_c(X, \mathcal{O}_E)$ \cite{ekedahl2007adic,laszlo2008six}. There is yet another construction of $D^b_c(X,\overline{\mathbb{Q}}_{\ell})$ coming from the pro-étale topology introduced in \cite{bhatt2015pro} which is the most natural setting for dg-purposes. In fact, in the pro-étale topology, the constructible derived category $D^b_c(X,\overline{\mathbb{Q}}_{\ell})$ does arise as a full triangulated subcategory of a genuine derived category and thus it is naturally dg-enhanced. For this reason, we choose to work with $D^b_c(X,\overline{\mathbb{Q}}_{\ell})$ in its pro-étale realization (which will be reviewed in \cref{Section: proetale topology}). To get a "purity implies formality" result, one also has to lift the Frobenius action on morphisms of sheaves to the dg-level as in \cite[Appendix A]{polishchuk2019semiorthogonal}. We explain a variant of this Frobenius lifting argument in \cref{Section: The Frobenius action on dg-algebras}. Once the Frobenius action is lifted to the dg-level, formality will be a consequence of an algebraic result in \cite[Theorem B.1.1]{polishchuk2019semiorthogonal}. While it would certainly be possible to directly apply the results from \cite{polishchuk2019semiorthogonal} in our setting without referring to the pro-étale topology, we hope that our alternative approach clarifies some of the technical difficulties.
\subsection*{Relation to other work}
Formality has been discussed in many settings of representation theory such as the Springer correspondence \cite{rider2013formality,polishchuk2019semiorthogonal,rider2021formality,
eberhardt2022motivic} or in the context of flag varieties and Koszul duality \cite{beilinson1996koszul,schnurer2011equivariant} where formality has also been studied for modular coefficients \cite{riche2014modular,achar2014modular}. It would also be interesting to prove similar formality results for graded Hecke algebras at central characters. These algebras can be used to study a wider range of representations of $p$-adic groups such as unipotent representations \cite{lusztig1995classification}. In terms of geometry, graded Hecke algebras arise as certain Ext-algebras in an \textit{equivariant} derived category of constructible sheaves \cite{lusztig1995cuspidal,aubert2018graded}. Some formality results in this direction can be found in \cite{solleveld2022graded}. There also is a coherent categorical Deligne-Langlands correspondence \cite{ben2020coherent} which works without fixing a central character but replaces constructible sheaves with a certain category of coherent sheaves (see also the conjectures in \cite{hellmann2020derived}). The relation between the constructible and the coherent side is discussed in \cite{ben2023between}.

\subsection*{Acknowledgments}
I am very grateful to my supervisor Kevin McGerty for his guidance and support and I would like to thank him for many enlightening discussions and conversations. I would also like to thank Dan Ciubotaru for many useful conversations and Ruben La and Emile Okada for many useful discussions on affine Hecke algebras and perverse sheaves.
\section{Springer geometry}
Let $\overline{\mathbb{F}}$ be an algebraically closed field and $\ell$ a prime number which is invertible in $\overline{\mathbb{F}}$. For any variety $X$ over $\overline{\mathbb{F}}$, we can consider the constructible derived category of $\overline{\mathbb{Q}}_{\ell}$-complexes $D^b_c(X, \overline{\mathbb{Q}}_{\ell} )$ as defined in \cite{deligne1980conjecture,beilinson2018faisceaux,bhatt2015pro}. The triangulated category $D^b_c(X, \overline{\mathbb{Q}}_{\ell} )$ comes with the usual six (derived) operations denoted by $f^*,f_*,f_!,f^!, \otimes^L$ and $\RIHom$. Moreover, there is a standard $t$-structure on $D^b_c(X, \overline{\mathbb{Q}}_{\ell} )$ with cohomology functor $H^*$ and heart $\Sh_c(X, \overline{\mathbb{Q}}_{\ell} )$. The structure map of $X$ will be denoted by $a: X \rightarrow \{\pt \}$. Let $\textbf{1}_X \in \Sh_c(X, \overline{\mathbb{Q}}_{\ell} )$ be the constant sheaf and $\omega_X:= a^! \textbf{1}_{\pt} \in D^b_c(X, \overline{\mathbb{Q}}_{\ell} )$ the dualizing complex. We denote by $\Perv(X) \subset D^b_c(X, \overline{\mathbb{Q}}_{\ell})$ the category of perverse sheaves on $X$ and by ${}^p H^*$ the perverse cohomology functor.
\subsection{Borel-Moore homology}
In this section we recall a few basic facts about Borel-Moore homology \cite{laumon1976homologie}. The $i$-th Borel-Moore homology is the $\overline{\mathbb{Q}}_{\ell}$-vector space
\begin{equation*}
H_i(X, \overline{\mathbb{Q}}_{\ell}) := \Hom^{-i}(\textbf{1}_X, \omega_X) = H^{-i} (a_*\omega_X).
\end{equation*}
It can be shown that $H_i(X, \overline{\mathbb{Q}}_{\ell})$ is concentrated in degrees $0,1,...,2 \dim X$. Moreover, we have the Künneth formula
\begin{equation*}
H_*(X \times Y , \overline{\mathbb{Q}}_{\ell} ) \cong  H_*(X , \overline{\mathbb{Q}}_{\ell}) \otimes H_*(Y , \overline{\mathbb{Q}}_{\ell}).
\end{equation*}
If $i: Y \hookrightarrow X$ is a closed immersion with complement $j: U \hookrightarrow X$ and $\mathcal{F} \in D^b_c(X, \overline{\mathbb{Q}}_{\ell})$, there is a canonical distinguished triangle
\begin{equation}\label{eq: open closed distinguised triangle}
i_* i^!\mathcal{F} \rightarrow \mathcal{F} \rightarrow j_*j^!\mathcal{F} \rightarrow [1]
\end{equation}
which is natural in $\mathcal{F}$. For $\mathcal{F} = \omega_X$ this induces a long exact sequence on Borel-Moore homology
\begin{equation}\label{eq: LES for Borel-Moore}
\cdots \rightarrow H_{i+1}(U, \overline{\mathbb{Q}}_{\ell}) \rightarrow H_i(Y, \overline{\mathbb{Q}}_{\ell}) \rightarrow H_i(X, \overline{\mathbb{Q}}_{\ell}) \rightarrow H_i(U, \overline{\mathbb{Q}}_{\ell}) \rightarrow H_{i-1}(Y, \overline{\mathbb{Q}}_{\ell}) \rightarrow \cdots .
\end{equation}
Let $p:\tilde{X} \rightarrow X$ be a smooth morphism of relative dimension $d$. The adjoint pair $(p^*, p_*)$ gives rise to a canonical morphism
\begin{equation*}
 \omega_X \rightarrow p_*p^*\omega_X  =  p_*p^!\omega_X [-2d]=  p_*\omega_{\tilde{X}}[-2d].
\end{equation*}
This induces a `smooth pullback' map on Borel-Moore homology
\begin{equation*}
 H_i(X, \overline{\mathbb{Q}}_{\ell}) \rightarrow H_{i+2d}(\tilde{X}, \overline{\mathbb{Q}}_{\ell}).
\end{equation*}
The naturality of the distinguished triangle in \eqref{eq: open closed distinguised triangle} implies that smooth pullback is compatible with the long exact sequence from \eqref{eq: LES for Borel-Moore}, i.e. if $Y \subset X$ is a closed subvariety with open complement $U$ and $\tilde{Y} := p^{-1}(Y)$, $\tilde{U} := p^{-1}(U)$, the following diagram commutes:
\begin{equation}\label{eq: smooth pullback is compatible with LES}
\begin{tikzcd}[sep=0.2in]
\cdots \arrow[r] &  H_{i+2d}(\tilde{Y}, \overline{\mathbb{Q}}_{\ell}) \arrow[r] & H_{i+2d}(\tilde{X}, \overline{\mathbb{Q}}_{\ell}) \arrow[r] & H_{i+2d}(\tilde{U}, \overline{\mathbb{Q}}_{\ell}) \arrow[r] & H_{i+2d-1}(\tilde{Y}, \overline{\mathbb{Q}}_{\ell}) \arrow[r] & \cdots \\
\cdots \arrow[r] &  H_i(Y, \overline{\mathbb{Q}}_{\ell}) \arrow[r] \arrow[u]    & H_i(X, \overline{\mathbb{Q}}_{\ell}) \arrow[r] \arrow[u]    & H_i(U, \overline{\mathbb{Q}}_{\ell}) \arrow[r] \arrow[u]    & H_{i-1}(Y, \overline{\mathbb{Q}}_{\ell}) \arrow[r] \arrow[u]  & \cdots.
\end{tikzcd}  
\end{equation}
\begin{Lemma}\label{Lemma: affine fibrations induce isomorphisms on homology}
Let $p:\tilde{X} \rightarrow X$ be a Zariski locally trivial fibration with affine fiber $\mathbb{A}^d$. Then the smooth pullback map $ H_i(X, \overline{\mathbb{Q}}_{\ell}) \rightarrow H_{i+2d}(\tilde{X}, \overline{\mathbb{Q}}_{\ell})$ is an isomorphism for all $i \in \mathbb{Z}$.
\end{Lemma}
\begin{proof}
Using \eqref{eq: smooth pullback is compatible with LES} and the five lemma, one can reduce to the case where $p$ is trivial. By the Künneth formula, it suffices to consider the case where $p: \mathbb{A}^d \rightarrow \pt$. Note that $H_i(\mathbb{A}^d, \overline{\mathbb{Q}}_{\ell}) = 0$ for $i \neq 2d $ and $H_{2d}(\mathbb{A}^d, \overline{\mathbb{Q}}_{\ell}) = \overline{\mathbb{Q}}_{\ell}$. The claim now follows since the smooth pullback map $H_0(\pt, \overline{\mathbb{Q}}_l) \rightarrow H_{2d}(\mathbb{A}^d, \overline{\mathbb{Q}}_{\ell})  $ is non-zero.
\end{proof}
If $X$ is smooth and connected, the dualizing complex is given by
\begin{equation*}
\omega_X= \textbf{1}_X [2 \dim X].
\end{equation*}
The fundamental class of $X$ is the distinguished element
\begin{equation*}
\id_{\textbf{1}_X} \in \Hom^0(\textbf{1}_X, \textbf{1}_X) = \Hom^{-2 \dim X}(\textbf{1}_X, \omega_X) = H_{2 \dim X} (X, \overline{\mathbb{Q}}_{\ell})
\end{equation*}
also denoted by $[X] \in H_{2 \dim X}(X, \overline{\mathbb{Q}}_{\ell})$. More generally, if $X$ is irreducible one can define the fundamental class as follows: Pick a smooth open subset $U \subset X$. Then the long exact sequence \eqref{eq: LES for Borel-Moore} induces an isomorphism $H_{2\dim X}(U, \overline{\mathbb{Q}}_{\ell}) \cong H_{2 \dim X}(X, \overline{\mathbb{Q}}_{\ell}) $. The image of $[U]$ under this isomorphism defines a distinguished element $[X] \in H_{2 \dim X}(X, \overline{\mathbb{Q}}_{\ell})$ called the fundamental class of $X$. It can be shown that $[X]$ does not depend on the choice of $U$. If $Y \subset X$ is an irreducible closed subvariety, the image of the fundamental class of $Y$ under $H_{2 \dim Y}(Y, \overline{\mathbb{Q}}_{\ell}) \rightarrow H_{2 \dim Y}(X , \overline{\mathbb{Q}}_{\ell})$ also defines a fundamental class $[Y] \in H_{2 \dim Y}(X, \overline{\mathbb{Q}}_{\ell})$. If elements of this form span the vector space $H_i(X, \overline{\mathbb{Q}}_{\ell})$ for each $i \in \mathbb{Z}$, we say that $H_*(X, \overline{\mathbb{Q}}_{\ell})$ is \textit{spanned by fundamental classes}. Note that being spanned by fundamental classes implies that $H_i(X, \overline{\mathbb{Q}}_{\ell})= 0$ for $i$ odd. Let $Z_i(X)$ be the free abelian group on the set of $i$-dimensional irreducible closed subvarieties of $X$ and let $A_i(X) = Z_i(X)/ \sim_{Rat}$ be the Chow group (c.f. \cite{fulton2013intersection}). The fundamental class construction gives rise to a cycle class map
\begin{align*}
cl_X: Z_i(X) & \rightarrow H_{2i}(X, \overline{\mathbb{Q}}_{\ell}) \\
		[Y] & \mapsto [Y]
\end{align*}
which descends to the Chow group
\begin{equation}\label{eq: cycle class map}
cl_X: A_i(X) \rightarrow H_{2i}(X, \overline{\mathbb{Q}}_{\ell})
\end{equation} 
by \cite[Théorème 6.3]{laumon1976homologie}. The open-closed exact sequence and smooth pullback map from Borel-Moore homology have analogues for Chow groups: For any $Y \subset X$ closed with complement $U = X\backslash Y$ there is an exact sequence
\begin{equation*}
A_i(Y) \rightarrow A_i(X) \rightarrow A_i(U) \rightarrow 0
\end{equation*}
and for $\tilde{p} : \tilde{X} \rightarrow X$ smooth (or more generally flat) of relative dimension $d$, there is a pullback map
\begin{equation*}
 A_i(X) \rightarrow A_{i+d}(X).
\end{equation*}
The cycle class map $cl_X$ is functorial with respect to these constructions \cite[Théoreme 6.1]{laumon1976homologie}. We define
\begin{equation*}
A_i(X)_{\overline{\mathbb{Q}}_{\ell}} := A_i(X) \otimes_{\mathbb{Z}} \overline{\mathbb{Q}}_{\ell}.
\end{equation*}
Following \cite[1.7]{de1988homology}, we say that a variety $X$ has property (S) if
\begin{itemize}
\item $H_i(X , \overline{\mathbb{Q}}_{\ell}) = 0$ for $i$ odd;
\item $cl_X: A_i(X)_{\overline{\mathbb{Q}}_{\ell}}   \rightarrow H_{2i}(X, \overline{\mathbb{Q}}_{\ell})$ is an isomorphism for all $i \in \mathbb{Z}$.
\end{itemize}
Note that if $X$ has property (S), then $H_*(X, \overline{\mathbb{Q}}_{\ell})$ is spanned by fundamental classes. The following two lemmas are $\overline{\mathbb{Q}}_{\ell}$-versions of \cite[Lemma 1.8, 1.9]{de1988homology}.
\begin{Lemma}\label{Lemma: Spanned by fundamental classes compatible with long exact sequence}
Let $Y \subset X$ be a closed subvariety with complement $U = X \backslash Y$. If $U$ and $Y$ have property (S) then $X$ also has property (S).
\end{Lemma}
\begin{proof}
The groups $H_i(Y, \overline{\mathbb{Q}}_{\ell})$ and $H_i(U, \overline{\mathbb{Q}}_{\ell})$ vanish for $i$ odd. Using the long exact sequence \eqref{eq: LES for Borel-Moore}, we deduce that $H_i(X, \overline{\mathbb{Q}}_{\ell})=0$ for $i$ odd. Moreover, for any $i \in \mathbb{Z}$, we get a commutative diagram with exact rows
\begin{equation*}
\begin{tikzcd}
            & A_i(Y)_{\overline{\mathbb{Q}}_{\ell}}  \arrow[r] \arrow[d, "\wr"] & A_i(X)_{\overline{\mathbb{Q}}_{\ell}}  \arrow[r] \arrow[d] & A_i(U)_{\overline{\mathbb{Q}}_{\ell}}  \arrow[r] \arrow[d, "\wr"] & 0 \\
0 \arrow[r] & {H_{2i}(Y, \overline{\mathbb{Q}}_{\ell})} \arrow[r]               & {H_{2i}(X, \overline{\mathbb{Q}}_{\ell})} \arrow[r]        & {H_{2i}(U, \overline{\mathbb{Q}}_{\ell})} \arrow[r]               & 0.
\end{tikzcd}
\end{equation*}
By the five lemma, this implies that the map $A_i(X)_{\overline{\mathbb{Q}}_{\ell}}  \rightarrow {H_{2i}(X, \overline{\mathbb{Q}}_{\ell})}$ is an isomorphism.
\end{proof}
\begin{Lemma}\label{Lemma: affine fibrations preserve spanned by fundamental class property}
Let $p:\tilde{X} \rightarrow X$ be a Zariski locally trivial fibration with fiber $\mathbb{A}^d$. If $X$ has property (S) then $\tilde{X}$ also has property (S).
\end{Lemma}
\begin{proof}
By \cref{Lemma: affine fibrations induce isomorphisms on homology}, the pullback map $H_i(X, \overline{\mathbb{Q}}_{\ell}) \rightarrow H_{i+2d}(\tilde{X}, \overline{\mathbb{Q}}_{\ell})$ is an isomorphism. In particular, $H_i(X, \overline{\mathbb{Q}}_{\ell})  = 0$ for $i$ odd. Moreover, for any $i \in \mathbb{Z}$ we get a commutative diagram
\begin{equation*}
\begin{tikzcd}
A_{i+d}(\tilde{X})_{\overline{\mathbb{Q}}_{\ell}} \arrow[r]      & {H_{2i+2d}(\tilde{X}, \overline{\mathbb{Q}}_{\ell})}    \\
A_i(X)_{\overline{\mathbb{Q}}_{\ell}} \arrow[r, "\sim"] \arrow[u] & {H_{2i}(X,\overline{\mathbb{Q}}_{\ell})}. \arrow[u, "\wr"]
\end{tikzcd}
\end{equation*}
As a consequence, the map $A_i(X)_{\overline{\mathbb{Q}}_{\ell}} \rightarrow A_{i+d}(\tilde{X})_{\overline{\mathbb{Q}}_{\ell}} $ is injective. It is also surjective by \cite[Proposition 1.9]{fulton2013intersection} an thus an isomorphism. This implies that $A_{i+d}(\tilde{X})_{\overline{\mathbb{Q}}_{\ell}} \rightarrow H_{2i+2d}(\tilde{X}, \overline{\mathbb{Q}}_{\ell})$ is also an isomorphism.
\end{proof}
\subsection{Morphisms of Springer type}\label{Section: morphisms of Springer type}
We now introduce a general setting of `Springer geometry' which we want to study in this paper. Let $G$ be a connected reductive group defined over $\overline{\mathbb{F}}$. We fix a maximal torus and a Borel subgroup $T \subset B \subset G$ with Weyl group $W = N_G(T)/T$. The corresponding flag variety will be denoted by $\mathcal{B} = G/B$. Given a $G$-representation $V$ and a finite collection $\{V^i \subset V \mid i \in I \}$ of $B$-stable subspaces, we define for each $i\in I$ the $G$-variety
\begin{equation*}
\tilde{V}^i := G \times^B V^i.
\end{equation*}
This comes with two $G$-equivariant morphisms
\begin{equation*}
\begin{tikzcd}
& \tilde{V}^i \arrow[ld, "\pi_i"'] \arrow[rd, "\mu_i"] & \\
\mathcal{B} & & V
\end{tikzcd}
\end{equation*}
where $\mu_i(g,v) = g \cdot v$ and $\pi_i(g,v) = gB$.
\begin{Definition}\label{Definition: morphism of Springer type}
We call $\mu_i$ a \textit{morphism of Springer type}.
\end{Definition}
\begin{Example}
Let $V = \mathfrak{g}$ be the Lie algebra of $G$, $\mathfrak{b}$ the Lie algebra of $B$ and $\mathfrak{n} = [\mathfrak{b}, \mathfrak{b}]$. Then the morphism of Springer type corresponding to the $B$-stable subspace $\mathfrak{n} \subset \mathfrak{g}$ is the Springer resolution
\begin{equation*}
\tilde{\mathcal{N}} = G \times^B \mathfrak{n} \rightarrow \mathcal{N}\subset \mathfrak{g}.
\end{equation*}
For $\mathfrak{b} \subset \mathfrak{g}$ we recover the Grothendieck-Springer alteration $G\times^B \mathfrak{b} \rightarrow \mathfrak{g}$. Other important examples of morphisms of Springer type show up in the representation theory of affine Hecke algebras (see \cite{kazhdan1987proof,chriss2009representation} and \cref{Corollary: the sqSpringer resolution is of Springer type}).
\end{Example}
For $i,j \in I$, we consider the Steinberg variety
\begin{equation*}
Z^{ij} := \tilde{V}^i \times_V \tilde{V}^j.
\end{equation*}
This comes with the projection map $\pi_i \times \pi_j : Z^{ij} \rightarrow \mathcal{B} \times \mathcal{B}$. Consider the orbit partition
\begin{equation*}
\mathcal{B} \times \mathcal{B} = \bigsqcup_{w \in W} Y_w
\end{equation*}
with
\begin{equation*}
Y_w := G\cdot (eB, \dot{w}B) \cong G/(B\cap \dot{w}B\dot{w}^{-1})
\end{equation*}
where $\dot{w}\in N_G(T)$ is a lift of $w \in W$. This induces a partition of $Z^{ij}$ into locally closed subvarieties
\begin{equation*}
Z^{ij} = \bigsqcup_{w\in W} Z^{ij}_w
\end{equation*}
where
\begin{equation*}
Z^{ij}_w:=(\pi_i \times \pi_j)^{-1}(Y_w).
\end{equation*}

\begin{Lemma}\label{Lemma: geometric properties of morphisms of Springer type}
\begin{enumerate}
\item The morphism $\mu_i : \tilde{V}^i \rightarrow V$ is proper;
\item The morphism $\pi_i : \tilde{V}^i \rightarrow \mathcal{B}$ is a Zariski vector bundle with fiber $V^i$;
\item For each $w \in W$ the morphism $\pi_i \times \pi_j : Z^{ij}_w \rightarrow Y_w$ is a Zariski vector bundle with fiber $V^i \cap \dot{w} V^j$;
\item The first projection $p_1:Y_w \rightarrow \mathcal{B}$ is a Zariski locally trivial fibration with fiber $\mathbb{A}^{l(w)}$.
\end{enumerate}
\end{Lemma}
\begin{proof}
The map $\mu_i$ can be factored into a closed immersion followed by a projection:
\begin{align*}
\tilde{V}^i & \hookrightarrow \mathcal{B} \times V \overset{p_2}{\longrightarrow} V\\
(g,v) & \mapsto (gB, gv).
\end{align*}
Since $\mathcal{B}$ is projective, this implies that $\mu_i$ is proper. The local triviality in (2)-(4) follows from standard results about quotients (c.f. \cite[§I.5.16]{jantzen2003representations}). The respective fibers are easily computed.
\end{proof}
\begin{Proposition}\label{Prop: H_*(Z) is spanned by fundamental classes}
The variety $Z^{ij}$ has property (S). In particular, $H_*(Z^{ij}, \overline{\mathbb{Q}}_{\ell})$ is spanned by fundamental classes.
\end{Proposition}
\begin{proof}
By \cref{Lemma: geometric properties of morphisms of Springer type} the maps $Z^{ij}_w \rightarrow Y_w \rightarrow \mathcal{B}$ are locally trivial fibrations with affine fibers. It is well known that $\mathcal{B}$ has property (S) (in fact, this follows from \cref{Lemma: Spanned by fundamental classes compatible with long exact sequence} and the decomposition of $\mathcal{B}$ into Bruhat cells). By \cref{Lemma: affine fibrations preserve spanned by fundamental class property} this implies that $Z^{ij}_w$ has property (S) for all $w \in W$. Pick a total order $\le$ on $W$ extending the Bruhat order and define
\begin{align*}
Z^{ij}_{\le w} &:= \bigsqcup_{y \le w} Z^{ij}_y \\
Z^{ij}_{< w} &:= \bigsqcup_{y < w} Z^{ij}_y.
\end{align*}
Note that $Z^{ij}_{< w} = Z^{ij}_{\le w'}$ where $w' \in W$ is the maximal element with $w' < w$. We show by induction along the total order on $W$ that $Z^{ij}_{\le w}$ has property (S). We have already proved the claim for $Z^{ij}_{\le e} = Z^{ij}_e$. Now assume we have shown the claim for each $y < w$. Note that $Z^{ij}_{< w} \subset Z^{ij}_{\le w}$ is closed with open complement $Z^{ij}_{w}$. The varieties $Z^{ij}_{< w}$ and $Z^{ij}_w$ have property (S), so $Z^{ij}_{\le w}$ also has property (S) by \cref{Lemma: Spanned by fundamental classes compatible with long exact sequence}. This completes the induction. Note that $Z^{ij} = Z^{ij}_{\le w_0}$ where $w_0 \in W$ is the longest element. Hence, $Z^{ij}$ has property (S).
\end{proof}
Since $\tilde{V}^i$ is smooth, there is the constant perverse sheaf
\begin{equation*}
\mathcal{C}_{\tilde{V}^i} := \textbf{1}_{\tilde{V}^i} [\dim \tilde{V}^i] \in \Perv(\tilde{V}^i).
\end{equation*}
We define the \textit{Springer sheaves}
\begin{equation*}
\Sbf^i := (\mu_i)_* \mathcal{C}_{\tilde{V}^i} \in D^b_c(V, \overline{\mathbb{Q}}_{\ell})
\end{equation*}
and
\begin{equation*}
\Sbf := \bigoplus_{i\in I} \Sbf^i.
\end{equation*}
The morphism $\mu_i$ is proper by \cref{Lemma: geometric properties of morphisms of Springer type}. Hence, the decomposition theorem \cite{beilinson2018faisceaux} implies that $\Sbf^i$ (and thus also $\Sbf$) is a semisimple complex. In other words, we have
\begin{equation*}
\Sbf = \bigoplus_{j=1}^n X_j [k_j]
\end{equation*}
for some simple perverse sheaves $X_j \in \Perv(V)$ and integers $k_j \in \mathbb{Z}$. We define the \textit{Springer category}
\begin{equation}\label{eq: definition of springer category}
D_{\Spr}(V, \overline{\mathbb{Q}}_{\ell}):= \langle X_1, ...,X_n \rangle_{\Delta}
\end{equation}
to be the smallest full triangulated subcategory of $D^b_c(V, \overline{\mathbb{Q}}_{\ell})$ that is closed under isomorphisms and contains the simple perverse sheaves $X_1,...,X_n$. Our main goal is to give an algebraic description of the Springer category.
\section{Purity}
In this section, we show that the canonical Frobenius action on the space of morphisms between any two Springer sheaves is pure.
\subsection{The Frobenius action}\label{Section: Frobenius action}
Let $X_0$ be a variety defined over a finite field $\mathbb{F}_q$ with structure map $a: X_0 \rightarrow \Spec(\mathbb{F}_q)$ and let
\begin{equation*}
X:= X_0 \otimes_{\mathbb{F}_q} \overline{\mathbb{F}}_q.
\end{equation*}
Denote by $D^b_c(X_0, \overline{\mathbb{Q}}_{\ell})$ the constructible derived category of $\overline{\mathbb{Q}}_{\ell}$-complexes on $X_0$. For any $\mathcal{F}_0 \in D^b_c(X_0, \overline{\mathbb{Q}}_{\ell})$ the pullback of $\mathcal{F}_0$ along $X \rightarrow X_0$ will be denoted by $\mathcal{F} \in D^b_c(X, \overline{\mathbb{Q}}_{\ell})$. We also define
\begin{equation}\label{eq: def of uHom on constructible sheaves}
\uHom^i(\mathcal{F}_0, \mathcal{G}_0) := H^i(a_*\RIHom(\mathcal{F}_0, \mathcal{G}_0)) \in \Sh_c(\Spec(\mathbb{F}_q), \overline{\mathbb{Q}}_{\ell})
\end{equation}
for any $\mathcal{F}_0 ,\mathcal{G}_0 \in D^b_c(X_0, \overline{\mathbb{Q}}_{\ell})$. The category $\Sh_c(\Spec(\mathbb{F}_q), \overline{\mathbb{Q}}_{\ell})$ is equivalent to the category of finite-dimensional continuous $\Gal(\overline{\mathbb{F}}_q/ \mathbb{F}_q)$-representations. Hence, we can consider $\uHom^i(\mathcal{F}_0, \mathcal{G}_0)$ as a $\overline{\mathbb{Q}}_{\ell}$-vector space equipped with a canonical Frobenius action (see \cref{Section: The Frobenius action on dg-algebras} for a construction of this action). Forgetting this action recovers the vector space $\Hom^i(\mathcal{F}, \mathcal{G})$. Note that the sheaves $\textbf{1}_X, \omega_X \in D^b_c(X, \overline{\mathbb{Q}}_{\ell})$ have canonical $\mathbb{F}_q$-versions $\textbf{1}_{X_0}$ and $\omega_{X_0} = a^! \textbf{1}_{\Spec(\mathbb{F}_q)}$. This induces a canonical Frobenius action on $H_*(X, \overline{\mathbb{Q}}_{\ell})$ coming from the Frobenius action on
\begin{equation*}
\underline{H}_i(X, \overline{\mathbb{Q}}_{\ell}) := \uHom^{-i}(\textbf{1}_{X_0}, \omega_{X_0}) = H^i(a_*\omega_{X_0}).
\end{equation*}
It can be shown that the constructions on Borel-Moore homology from the previous section are compatible with the Frobenius action. For example, if $i: Y_0 \hookrightarrow X_0$ is a closed immersion with open complement $j: U_0 \hookrightarrow X_0$, there is a distinguished triangle
\begin{equation*}
i_*i^! \omega_{X_0} \rightarrow \omega_{X_0}  \rightarrow j_*j^! \omega_{X_0}  \rightarrow [1]
\end{equation*}
inducing a long exact sequence
\begin{equation*}
\cdots \rightarrow \underline{H}_{i+1}(U, \overline{\mathbb{Q}}_{\ell}) \rightarrow \underline{H}_i(Y, \overline{\mathbb{Q}}_{\ell}) \rightarrow \underline{H}_i(X, \overline{\mathbb{Q}}_{\ell}) \rightarrow \underline{H}_i(U, \overline{\mathbb{Q}}_{\ell}) \rightarrow \underline{H}_{i-1}(Y, \overline{\mathbb{Q}}_{\ell}) \rightarrow \cdots
\end{equation*}
where all maps commute with the Frobenius action. For any $n \in \mathbb{Z}$, let
\begin{equation*}
\textbf{1}_{\Spec(\mathbb{F}_q)}(n) \in \Sh_c(\Spec(\mathbb{F}_q), \overline{\mathbb{Q}}_{\ell})
\end{equation*}
be the $n$-th Tate twist. This corresponds to the $1$-dimensional $\overline{\mathbb{Q}}_{\ell}$-vector space on which the (geometric) Frobenius element acts by multiplication with $q^{-n}$. Moreover, we define
\begin{align*}
\textbf{1}_{X_0}(n) & := a^* \textbf{1}_{\Spec(\mathbb{F}_q)}(n) \in \Sh_c(X_0, \overline{\mathbb{Q}}_{\ell}) \\
\mathcal{F}_0(n) &:= \mathcal{F}_0 \otimes^L \textbf{1}_{X_0}(n) \in D^b_c(X_0 , \overline{\mathbb{Q}}_{\ell})
\end{align*}
for any $\mathcal{F}_0 \in D^b_c(X_0 , \overline{\mathbb{Q}}_{\ell})$ and $n \in \mathbb{Z}$. If $X_0$ is smooth and connected, the dualizing complex on $X_0$ is given by
\begin{equation*}
\omega_{X_0} = \textbf{1}_{X_0}[2 \dim X](\dim X).
\end{equation*}
\begin{Lemma}\label{Lemma: Frobenius action on fundamental class}
Let $X$ be irreducible. Then Frobenius acts on $[X] \in \underline{H}_{2 \dim X}(X, \overline{\mathbb{Q}}_{\ell})$ by multiplication with $q^{- \dim X}$.
\end{Lemma}
\begin{proof}
Let $U_0 \subset X_0$ be a smooth open subset. Then we get a canonical element
\begin{align*}
\id_{\mathbf{1}_{U_0}} & \in \underline{\Hom}^0(\mathbf{1}_{U_0},\mathbf{1}_{U_0}) \\
& \cong \underline{\Hom}^0(\mathbf{1}_{U_0},\omega_{U_0}[-2\dim U](-\dim U)) \\
& = \underline{\Hom}^{-2\dim U}(\mathbf{1}_{U_0},\omega_{U_0})(-\dim U)
\end{align*}
which is Frobenius invariant and corresponds to the fundamental class $[U] \in H_{2 \dim U}(U, \overline{\mathbb{Q}}_{\ell})$ after forgetting the Frobenius action. Note that the Frobenius action on $ \underline{H}_{2 \dim U}(U, \overline{\mathbb{Q}}_{\ell})$ comes from the Frobenius action on $\underline{\Hom}^{-2\dim U}(\mathbf{1}_{U_0},\omega_{U_0})$. Thus, taking into account the Tate twist, we get that Frobenius acts on $ [U]$ by multiplication with $q^{- \dim U}$. The fundamental class $[X]$ is the inverse image of $[U]$ under the Frobenius equivariant restriction isomorphism $\underline{H}_{2 \dim X}(X, \overline{\mathbb{Q}}_{\ell}) \overset{\sim}{\rightarrow} \underline{H}_{2 \dim X}(U, \overline{\mathbb{Q}}_{\ell})$. Hence, Frobenius also acts on $[X]$ by multiplication with $q^{- \dim U} = q^{- \dim X}$.
\end{proof}
\begin{Corollary}\label{Cor: spanned by fundamental classes implies purity}
Assume that $H_*(X, \overline{\mathbb{Q}}_{\ell})$ is spanned by fundamental classes. Then for $\mathbb{F}_q$ large enough, Frobenius acts on $\underline{H}_i(X, \overline{\mathbb{Q}}_{\ell})$ by multiplication with $q^{- \frac{i}{2}}$.
\end{Corollary}
\begin{proof}
Let $Y_1,...,Y_n$ be irreducible closed subvarieties of $X$ such that the fundamental classes $[Y_1],...,[Y_n]$ span $H_*(X, \overline{\mathbb{Q}}_{\ell})$. For $\mathbb{F}_q$ large enough, we may assume that each of the $Y_i$ can be defined over $\mathbb{F}_q$. Then \cref{Lemma: Frobenius action on fundamental class} implies that Frobenius acts on $[Y_i] \in \underline{H}_{2 \dim Y_i}(X, \overline{\mathbb{Q}}_{\ell})$ by multiplication with $q^{- \dim Y_i}$. Since the $[Y_i]$ span $H_*(X, \overline{\mathbb{Q}}_{\ell})$, this proves the claim.
\end{proof}
\subsection{Frobenius and the Springer sheaf}
Let $\mu_i : \tilde{V}^i \rightarrow V$ ($i \in I$) be a finite collection of morphisms of Springer type (\cref{Definition: morphism of Springer type}) defined over $\overline{\mathbb{F}}_q$. Assume that $\mathbb{F}_q$ is large enough so that there are $\mathbb{F}_q$-forms $\mu_i: \tilde{V}^i_0 \rightarrow V_0$ for each of the $\mu_i$ ($i \in I$). Then $Z^{ij}$ can also be defined over $\mathbb{F}_q$ by considering $Z^{ij}_0 := \tilde{V}^i_0 \times_{V_0} \tilde{V}^j_0$.
\begin{Lemma}\label{Lemma: Frobenius action on H(Z)}
For $\mathbb{F}_q$ large enough, Frobenius acts on $\underline{H}_k(Z^{ij}, \overline{\mathbb{Q}}_{\ell})$ by multiplication with $q^{-\frac{k}{2}}$.
\end{Lemma}
\begin{proof}
This follows from \cref{Prop: H_*(Z) is spanned by fundamental classes} and \cref{Cor: spanned by fundamental classes implies purity}.
\end{proof}
We fix a square root $q^{\tfrac{1}{2}}$ of $q$ in $\overline{\mathbb{Q}}_l$. This corresponds to fixing a square root $\textbf{1}_{\Spec(\mathbb{F}_q)} (\tfrac{1}{2})$ of the Tate sheaf $\textbf{1}_{\Spec(\mathbb{F}_q)} (1)$. Then we can form the half integer Tate twists $\mathcal{F}_0(\tfrac{n}{2})$ for any $\mathcal{F}_0 \in D_c^b(X_0, \overline{\mathbb{Q}}_{\ell})$ and $n\in \mathbb{Z}$. Consider the constant perverse sheaf of weight $0$ on $\tilde{V}^i_0$
\begin{equation*}
\mathcal{C}_{\tilde{V}^i_0} := \mathbf{1}_{\tilde{V}^i_0} [d_i](\tfrac{d_i}{2}) \in D^b_c(\tilde{V}_0, \overline{\mathbb{Q}}_{\ell})
\end{equation*}
where $d_i = \dim \tilde{V}^i$ and let
\begin{equation*}
\mathbf{S}^i_0 := (\mu_i)_*\mathcal{C}_{\tilde{V}^i_0} \in D^b_c(V_0, \overline{\mathbb{Q}}_{\ell})
\end{equation*}
be the corresponding $\mathbb{F}_q$-Springer sheaf. We will need the following Frobenius equivariant version of \cite[Lemma 8.6.1]{chriss2009representation}.
\begin{Lemma}\label{Lemma: Frobenius equivariant iso between Ext and homology}
There is a (Frobenius equivariant) isomorphism
\begin{equation*}
\underline{\Hom}^k(\mathbf{S}^i_0,\mathbf{S}^j_0) \cong \underline{H}_{d_i+d_j-k}(Z^{ij}, \overline{\mathbb{Q}}_{\ell})(\tfrac{-d_i-d_j}{2}).
\end{equation*}
\end{Lemma}
\begin{proof}
Note that $\mu_i$ is proper over $\mathbb{F}_q$ since it is proper over $\overline{\mathbb{F}}_q$ (properness can be checked fpqc locally). Hence, we have $\mathbf{S}^i_0 = (\mu_i)_*\mathcal{C}_{\tilde{V}^i_0}= (\mu_i)_!\mathcal{C}_{\tilde{V}^i_0}$. Using base change with respect to the cartesian diagram
\begin{equation*}
\begin{tikzcd}
Z^{ij}_0 \arrow[d, "p"'] \arrow[r, "q"] & \tilde{V}^j_0 \arrow[d, "\mu_j"'] \\
\tilde{V}^i_0 \arrow[r, "\mu_i"]           & V_0                          
\end{tikzcd}
\end{equation*}
we get
\begin{align*}
\RIHom (\mathbf{S}^i_0,\mathbf{S}^j_0) & \cong \RIHom((\mu_i)_!\mathbf{1}_{ \tilde{V}^i_0},  (\mu_j)_* \mathbf{1}_{\tilde{V}^j_0})[d_j-d_i](\tfrac{d_j-d_i}{2})\\
&  \cong (\mu_i)_* \RIHom(\mathbf{1}_{\tilde{V}^i_0}, \mu_i^! (\mu_j)_* \mathbf{1}_{\tilde{V}^j_0})[d_j-d_i](\tfrac{d_j-d_i}{2})\\
&  \cong (\mu_i)_* \RIHom(\mathbf{1}_{\tilde{V}^i_0}, p_*q^! \mathbf{1}_{\tilde{V}^j_0})[d_j-d_i](\tfrac{d_j-d_i}{2}) \\
& \cong (\mu_i)_* p_* \RIHom(p^*\mathbf{1}_{\tilde{V}^i_0}, q^! \mathbf{1}_{\tilde{V}^j_0})[d_j-d_i](\tfrac{d_j-d_i}{2}) \\
& \cong (\mu_i)_* p_* \RIHom(p^*\mathbf{1}_{\tilde{V}^i_0}, q^! \omega_{\tilde{V}^j_0})[-d_j-d_i](\tfrac{-d_j-d_i}{2}) \\
& \cong (\mu_i)_* p_* \RIHom(\mathbf{1}_{Z^{ij}_0}, \omega_{Z^{ij}_0})[-d_j-d_i](\tfrac{-d_j-d_i}{2}).
\end{align*}
Applying $H^ka_*$ to this, we obtain
\begin{align*}
\underline{\Hom}^k(\mathbf{S}^i_0,\mathbf{S}^j_0) & = H^k(a_*\RIHom (\mathbf{S}^i_0,\mathbf{S}^j_0)) \\
& \cong H^{k-d_i-d_j} (a_* (\mu_i)_* p_* \RIHom(\mathbf{1}_{Z^{ij}_0}, \omega_{Z^{ij}_0}))(\tfrac{-d_i-d_j}{2}) \\
& \cong \underline{H}_{d_i+d_j-k}(Z^{ij}, \overline{\mathbb{Q}}_{\ell})(\tfrac{-d_i-d_j}{2}).
\end{align*}
\end{proof}
\begin{Corollary}\label{Cor: Purity of the Ext algebra}
For $\mathbb{F}_q$ large enough, Frobenius acts on $\underline{\Hom}^k(\mathbf{S}^i_0, \mathbf{S}^j_0)$ by multiplication with $q^{\frac{k}{2}}$.
\end{Corollary}
\begin{proof}
By \cref{Lemma: Frobenius action on H(Z)} Frobenius acts on $\underline{H}_{d_i+d_j-k}(Z^{ij})(\tfrac{-d_i-d_j}{2})$ by multiplication with
\begin{equation*}
q^{-\tfrac{d_i+d_j-k}{2}} \cdot q^{-\tfrac{-d_i-d_j}{2}} = q^{\tfrac{k}{2}}.
\end{equation*}
The claim now follows from \cref{Lemma: Frobenius equivariant iso between Ext and homology}.
\end{proof}
\section{Formality}
\subsection{Idempotent complete triangulated categories}
Recall that an additive category $\mathcal{A}$ is called \text{idempotent complete} if any idempotent $e : X \rightarrow X$ in $ \mathcal{T}$ splits. This is equivalent to the property that all idempotents have kernels (or cokernels). In an idempotent complete category, every idempotent $e: X \rightarrow X$ gives rise to a canonical decomposition
\begin{equation*}
X = \ker(e) \oplus \im(e).
\end{equation*}
If $\mathcal{A}' \subset \mathcal{A}$ is an additive subcategory that is closed under direct summands and $\mathcal{A}$ is idempotent complete, then $\mathcal{A}'$ is also idempotent complete. Moreover, for any additive category $\mathcal{A}$ one can define its \textit{idempotent completion} $\tilde{\mathcal{A}}$. This is an idempotent complete additive category with a fully faithful additive functor $\iota: \mathcal{A} \rightarrow \tilde{\mathcal{A}}$ and the following universal property: Any additive functor $\mathcal{A} \rightarrow \mathcal{C}$ with $\mathcal{C}$ idempotent complete factors uniquely (up to isomorphism) through $\iota$. For more details about idempotent completeness we refer to \cite[I.6]{karoubi1978k}.
\medskip

A triangulated category is called idempotent complete if its underlying additive category is idempotent complete. All triangulated categories we will encounter are idempotent complete thanks to the following well-known result (c.f. \cite[Lemma 1.6.8]{neeman2001triangulated} and \cite[Lemma 2.4]{balmer2001idempotent}).
\begin{Lemma}\label{Lemma: countable coproducts implies idempotent complete}
Any triangulated category with countable coproducts is idempotent complete. Moreover, the bounded below derived category $D^+(\mathcal{A})$ of an abelian category $\mathcal{A}$ is always idempotent complete.
\end{Lemma}
Let $\mathcal{T}$ be a triangulated category. For any object $X \in \mathcal{T}$, we denote by 
\begin{equation*}
\langle X \rangle_{\mathcal{T}, \Delta} = \langle X \rangle_{ \Delta}
\end{equation*}
the smallest full triangulated subcategory of $\mathcal{T}$ that contains $X$ and is closed under isomorphisms. Similarly, we denote by
\begin{equation*}
\langle X \rangle_{\mathcal{T}, \Delta, \inplus} = \langle X \rangle_{\Delta, \inplus}
\end{equation*}
the smallest full triangulated subcategory of $\mathcal{T}$ that contains $X$ and is closed under isomorphisms and direct summands.
\begin{Lemma}\label{Lemma: idempotent completion of a subcategory of an idempotent complete category}
Let $\mathcal{T}$ be an idempotent complete triangulated category and $X \in \mathcal{T}$. Then $\langle X \rangle_{\Delta, \inplus}$ is the idempotent completion of $\langle X \rangle_{\Delta}$.
\end{Lemma}
\begin{proof}
By \cite[Theorem I.6.12]{karoubi1978k} the idempotent completion of $\langle X \rangle_{\Delta}$ can be described as the full additive subcategory $\mathcal{C} \subset \mathcal{T}$ consisting of those objects $Y \in \mathcal{T}$ which are isomorphic to a direct summand of an object in $\langle X \rangle_{\Delta}$. In particular, we have $\mathcal{C} \subset \langle X \rangle_{\Delta, \inplus}$. The cone of a morphism in $\mathcal{C}$ is a direct summand of the cone of a morphism in $\langle X \rangle_{\Delta}$. Hence, the category $\mathcal{C}$ is closed under cones and thus $\mathcal{C} = \langle X \rangle_{\Delta, \inplus}$.
\end{proof}
We will also need the following idempotent complete version of Beilinson's Lemma (see \cite{beilinson1987derived} and \cite[Lemma 6]{schnurer2011equivariant}) .
\begin{Lemma}\label{Lemma: idempotent complete Beilinson's Lemma}
Let $F: \mathcal{T} \rightarrow \mathcal{T}'$ be a triangulated functor between idempotent complete triangulated categories. Let $X \in \mathcal{T}$ such that $F$ induces an isomorphism
\begin{equation*}
\Hom_{\mathcal{T}}(X,X[i]) \overset{\sim}{\rightarrow} \Hom_{\mathcal{T}'}(F(X),F(X)[i])
\end{equation*}
for all $i \in \mathbb{Z}$. Then $F$ restricts to an equivalence of triangulated categories
\begin{equation*}
\langle X \rangle_{ \Delta, \inplus} \cong \langle F(X) \rangle_{ \Delta, \inplus}.
\end{equation*}
\end{Lemma}
\begin{proof}
By a standard dévissage argument, $F$ induces an equivalence $\langle X \rangle_{ \Delta} \cong \langle F(X) \rangle_{ \Delta}$. This extends to an equivalence of the respective idempotent completions and thus to an equivalence $\langle X \rangle_{ \Delta, \inplus} \cong \langle F(X) \rangle_{ \Delta, \inplus}$ by \cref{Lemma: idempotent completion of a subcategory of an idempotent complete category}.
\end{proof}

\subsection{Derived categories and dg-algebras}\label{Section: Derived categories and dg-algebras}
Let $k$ be a commutative ring. For any dg-algebra $R$ over $k$, we denote by $R-\dgMod$ the category of (left) dg-modules over $R$. The homotopy category of dg-modules will be denoted by $K(R-\dgMod)$ and the derived category by $D(R-\dgMod)$. We write
\begin{equation*}
D_{\perf}(R-\dgMod) := \langle R \rangle_{D(R-\dgMod),\Delta, \inplus }
\end{equation*}
for the perfect derived category. For more details about dg-algebras and dg-modules we refer to \cite{keller1994deriving,bernstein2006equivariant}.
\medskip

Let $\mathcal{A}$ be a $k$-linear abelian category. Then for any two chain complexes $X^{\bullet}, Y^{\bullet} \in C(\mathcal{A})$, we can consider the Hom-complex $\Hom_{dg}^{\bullet}(X^{\bullet},Y^{\bullet}) \in C(k)$. Explicitly, this complex is defined as
\begin{equation*}
\Hom_{dg}^i(X^{\bullet},Y^{\bullet}) := \prod_{l+m=i} \Hom_{\mathcal{A}}(X^{-l}, Y^m)
\end{equation*}
with differential $d_{\Hom^{\bullet}_{dg}(X^{\bullet},Y^{\bullet})} (f) = d_{Y^{\bullet}} \circ f - (-1)^{|f|} f \circ d_{X^{\bullet}}$. Taking cohomology of this complex recovers the morphism space in the homotopy category:
\begin{equation*}
H^i( \Hom_{dg}^{\bullet}(X^{\bullet},Y^{\bullet}) )= \Hom_{K(\mathcal{A})}(X^{\bullet}, Y^{\bullet}[i]).
\end{equation*}
Component-wise composition defines a map
\begin{equation*}
\Hom_{dg}^i(Y^{\bullet},Z^{\bullet}) \otimes   \Hom_{dg}^j(X^{\bullet},Y^{\bullet})\rightarrow \Hom_{dg}^{i+j}(X^{\bullet},Z^{\bullet})
\end{equation*}
for any $X^{\bullet},Y^{\bullet},Z^{\bullet} \in C(\mathcal{A})$. This equips $\Hom^{\bullet}_{dg}(Y^{\bullet},Y^{\bullet})$ with the structure of a dg-algebra and $\Hom^{\bullet}_{dg}(X^{\bullet},Y^{\bullet})$ with the structure of a (left) dg-module over $\Hom^{\bullet}_{dg}(Y^{\bullet},Y^{\bullet})$. Furthermore, this gives rise to a functor
\begin{equation}\label{eq: chain complex functor Homdg(-,Y)}
\begin{aligned}
\Hom_{dg}^{\bullet}(-,Y^{\bullet}): C(\mathcal{A})^{\op} & \rightarrow \Hom^{\bullet}_{dg}(Y^{\bullet},Y^{\bullet})-\dgMod\\
X^{\bullet} & \mapsto \Hom^{\bullet}_{dg}(X^{\bullet},Y^{\bullet})
\end{aligned}
\end{equation}
which descends to a triangulated functor on the respective homotopy categories
\begin{equation}\label{eq: homotopy functor Homdg(-,Y)}
\Hom_{dg}^{\bullet}(-,Y^{\bullet}): K(\mathcal{A})^{\op} \rightarrow K(\Hom^{\bullet}_{dg}(Y^{\bullet},Y^{\bullet})-\dgMod).
\end{equation}
If $Y^{\bullet}$ is a $\mathcal{K}$-injective complex (e.g. a bounded below complex of injectives), this descends further to a triangulated functor on the corresponding derived categories
\begin{equation}\label{eq: derived functor Homdg(-,Y)}
\Hom_{dg}^{\bullet}(-,Y^{\bullet}): D(\mathcal{A})^{\op} \rightarrow D(\Hom^{\bullet}_{dg}(Y^{\bullet},Y^{\bullet})-\dgMod).
\end{equation}

Assume now that the abelian category $\mathcal{A}$ has enough injectives. Then for each $X \in D^+(\mathcal{A})$ we can pick a complex of injectives $I^{\bullet}_X \in C^+(\mathcal{A})$ representing $X$. Moreover, we can consider the associated dg-algebra
\begin{equation}\label{eq: definition of dg algebra associed to X}
R_X := \Hom^{\bullet}_{dg}(I_X^{\bullet},I_X^{\bullet})
\end{equation}
and the corresponding functor
\begin{equation}\label{eq:functor from derived category to dg-modules}
\Hom^{\bullet}_{dg}(-,I_X^{\bullet}): D(\mathcal{A})^{\op} \rightarrow D(R_X- \dgMod).
\end{equation}
The following statement is an idempotent complete version of \cite[Proposition 7]{schnurer2011equivariant}.
\begin{Lemma}\label{Cor: derived cat generated by X equivalent to perfect dgMod cat}
The functor from \eqref{eq:functor from derived category to dg-modules} induces an equivalence of triangulated categories
\begin{equation*}
\langle X  \rangle_{D^+(\mathcal{A})^{\op}, \Delta, \inplus} \cong D_{\perf}(R_X-\dgMod)
\end{equation*}
which sends $X$ to the free dg-module $R_X$.
\end{Lemma}
\begin{proof}
By construction, the functor from \eqref{eq:functor from derived category to dg-modules} sends $X$ to $\Hom^{\bullet}_{dg}(I_X^{\bullet}, I_X^{\bullet}) = R_X$. Moreover, this functor induces an isomorphism
\begin{align*}
\Hom_{D^+(A)}(X,X[i]) & \cong \Hom_{K(\mathcal{A})}(I^{\bullet}_X, I^{\bullet}_X[i]) \\
& \cong H^i(R_X) \\
& \cong \Hom_{K(R_X-\dgMod)}(R_X, R_X[i]) \\
& \cong \Hom_{D(R_X-\dgMod)}(R_X,R_X[i])
\end{align*}
for all $i \in \mathbb{Z}$. Here we have used for the first (resp. last) isomorphism that $I^{\bullet}_X$ (resp $R_X$) is a $\mathcal{K}$-injective complex (resp. $\mathcal{K}$-projective dg-module). Note that the categories $ D^+(\mathcal{A})^{\op}$ and $D(R_X-\dgMod)$ are idempotent complete by \cref{Lemma: countable coproducts implies idempotent complete} (using that $D(R_X-\dgMod)$ has countable coproducts). By \cref{Lemma: idempotent complete Beilinson's Lemma} this implies that \eqref{eq:functor from derived category to dg-modules} restricts to an equivalence of triangulated categories
\begin{equation*}
\langle X  \rangle_{D^+(\mathcal{A})^{\op}, \Delta, \inplus} \cong \langle R_X  \rangle_{D(R_X-\dgMod), \Delta, \inplus} = D_{\perf}(R_X-\dgMod).
\end{equation*}
\end{proof}
The dg-algebra $R_{X}$ crucially depends on the choice of the injective resolution $I_X^{\bullet}$ which makes it difficult to compute $R_X$ explicitly. The cohomology of $R_X$ on the other hand can be described more concretely as
\begin{equation*}
H^i(R_X) = H^i(\Hom^{\bullet}_{dg}(I_X^{\bullet}, I_X^{\bullet})) = \Hom^i(X,X).
\end{equation*}
Hence, we would like to replace the dg-algebra $R_X$ in \cref{Cor: derived cat generated by X equivalent to perfect dgMod cat} with its cohomology. Recall that a dg-algebra $R$ is called \textit{formal} if there is a chain of quasi-isomorphisms of dg-algebras
\begin{equation*}
R \leftarrow A_1 \rightarrow A_2 \leftarrow ... \rightarrow ... \leftarrow A_n \rightarrow H^*(R).
\end{equation*}
Here we consider $H^*(R)$ as a dg-algebra with vanishing differential. Any quasi-isomorphism of dg-algebras $A \rightarrow B$ induces an equivalence $D_{\perf}(B) \overset{\sim}{\rightarrow} D_{\perf}(A)$ which identifies the free dg-module $B$ with the free dg-module $A$.  Hence, if $R$ is formal, there is an equivalence of triangulated categories
\begin{equation*}
D_{\perf}(R-\dgMod) \cong D_{\perf}(H^*(R)- \dgMod).
\end{equation*}
In particular, if $R_X$ is formal \cref{Cor: derived cat generated by X equivalent to perfect dgMod cat} induces an equivalence
\begin{equation}\label{eq: equivalence of categories when Rx is formal}
\langle X \rangle _{D^+(\mathcal{A})^{\op},\Delta, \inplus} \cong D_{\perf}(\Hom^*(X,X)-\dgMod).
\end{equation}

The following theorem provides a useful criterion for formality.
\begin{Theorem}\cite{polishchuk2019semiorthogonal}\label{Lemma: purity implies formality}
Let $R$ be a dg-algebra over an algebraically closed field $k$ and $q \in k^{\times}$ not a root of unity. If $R$ can be equipped with a dg-algebra automorphism $F: R\rightarrow R$ such that $H^i(F)$ acts on $H^i(R)$ by multiplication with $q^i$ then $R$ is formal.
\end{Theorem}
\begin{proof}
This is proved for $k = \mathbb{C}$ in \cite[Theorem B.1.1]{polishchuk2019semiorthogonal} (with a slightly different condition on the action of $H^i(F)$). The same proof works for general $k$: By \cite[Theorem B.1.2]{polishchuk2019semiorthogonal}, we may assume that $F$ acts locally finitely on $R$. Hence, there is a generalized eigenspace decomposition $R^i= \bigoplus_{\alpha \in k^{\times}} R^i_{\alpha}$ which satisfies $R^i_{\alpha} \cdot R^j_{\beta} \subset R^{i+j}_{\alpha \beta}$. Using this decomposition, we can define a dgg-algebra (i.e. a $\mathbb{Z}^2$-graded algebra with a differential $d$ homogeneous of degree $(1,0)$ satisfying the Leibniz rule) via $\tilde{R} := \bigoplus_{i,j \in \mathbb{Z}} \tilde{R}^i_{q^j}$. Note that the inclusion $\tilde{R} \hookrightarrow R$ is a quasi-isomorphism. Moreover, the cohomology of the dgg-algebra $\tilde{R}$ lives in degrees $\{ (i,i) \mid i \in \mathbb{Z} \}$. Such dgg-algebras are known to be formal (c.f. \cite[Proposition 4]{schnurer2011equivariant}).
\end{proof}
\subsection{The pro-étale topology}\label{Section: proetale topology}
In this section, we recall the definition of the constructible derived category in the pro-étale topology and a few of its basic properties from \cite{bhatt2015pro}.
\begin{Definition}
A morphism of schemes $f: U \rightarrow X$ is called weakly \'{e}tale if both $f$ and its diagonal $\Delta_f: U \rightarrow U\times_X U $ are flat morphisms. Let $X_{\proet}$ be the category of schemes $U$ weakly  \'{e}tale over $X$. This becomes a site by declaring a family of maps $\{U_i \rightarrow U \}$ to be a covering if it is a covering in the fpqc topology.
\end{Definition}
The property of being weakly étale is stable under composition and base change. Moreover, any morphism in $X_{\proet}$ is weakly étale. We denote the category of sheaves of abelian groups on $X_{\proet}$ by $\Ab(X_{\proet})$. Let $f:X \rightarrow Y$ be a morphism of schemes. Then $f$ induces a morphism of sites $X_{\proet} \rightarrow Y_{\proet}$ via the functor
\begin{align*}
Y_{\proet} &\rightarrow X_{\proet} \\
V & \mapsto V \times_Y X.
\end{align*}
There is a corresponding pair of adjoint functors
\begin{equation*}
\begin{tikzcd}
\Ab(Y_{\proet}) \arrow[r, "f^{-1}", shift left]  & \Ab(X_{\proet}) \arrow[l, "f_*", shift left].
\end{tikzcd}
\end{equation*}
Let $\Lambda$ be a topological ring. By \cite[Lemma 4.2.12]{bhatt2015pro} there is a sheaf of rings $\Lambda_X$ on $X_{\proet}$ defined by
\begin{equation*}
\Lambda_X(U) := \Hom_{\cont}(U,\Lambda).
\end{equation*}
If $\Lambda$ is totally disconnected and $U$ is quasi-compact, we get
\begin{equation}\label{eq: constant sheaf evaluated at quasi compact covers}
\Lambda_X(U) = \Hom_{\cont}(\pi_0(U),\Lambda)
\end{equation}
where $\pi_0(U)$ is the space of connected components. If $\Lambda$ is discrete, $\Lambda_X$ is the constant sheaf with values in $\Lambda$. For any morphism of schemes $f :X \rightarrow Y$, there is a canonical map $f^{\#}: \Lambda_Y \rightarrow f_*\Lambda_X$ and hence a morphism of ringed sites
\begin{equation*}
(X_{\proet}, \Lambda_X) \rightarrow (Y_{\proet}, \Lambda_Y).
\end{equation*}
Thus we get corresponding pairs of adjoint functors
\begin{equation*}
\begin{tikzcd}
\Sh(Y_{\proet}, \Lambda_Y) \arrow[r, "f^*", shift left]  & \Sh(X_{\proet}, \Lambda_X) \arrow[l, "f_*", shift left] \\
D(Y_{\proet}, \Lambda_Y) \arrow[r, "f^*", shift left]  & D(X_{\proet}, \Lambda_X) \arrow[l, "f_*", shift left]. 
\end{tikzcd}
\end{equation*}
Here we adhere to the derived convention, i.e. the functors $f^*$ and $f_*$ are understood to be derived when applied to complexes. From now on, we fix a prime number $\ell$ and make the following assumption on the coefficient ring $\Lambda$.
\begin{Assumption}\label{Assumption: coefficient ring Lambda}
The ring $\Lambda$ is of one of the following forms:
\begin{enumerate}
\item $\Lambda = \mathbb{F}$ is a finite field of characteristic $\ell$;
\item $\Lambda= \mathcal{O}_E$ is the ring of integers in an algebraic extension $E/ \mathbb{Q}_{\ell}$;
\item $\Lambda = E$ is an algebraic extension $E / \mathbb{Q}_{\ell}$.
\end{enumerate}
We equip $\Lambda$ with the $\ell$-adic topology.
\end{Assumption}
Note that the assumption above includes the ring $\Lambda = \overline{\mathbb{Q}}_{\ell}$. For any algebraic extension $E/\mathbb{Q}_{\ell}$, we denote by $\mathcal{O}_E$ the ring of integers in $E$. If $E/\mathbb{Q}_{\ell}$ is finite, we pick a uniformizer $\omega \in \mathcal{O}_E$.
\begin{Example}\label{Example: explicit description of sheaf of ring in pro-etale topology}
\cite[Lemma 6.8.2]{bhatt2015pro} gives an alternative description of $\Lambda_X$ under \cref{Assumption: coefficient ring Lambda}: If $\Lambda$ is a finite ring (e.g. $\Lambda = \mathbb{F}$ a finite field) then $\Lambda_X$ is just the constant sheaf on $X_{\proet}$ with values in $\Lambda$. In the remaining cases, $\Lambda_X$ can be constructed as follows:
\begin{align*}
\mathcal{O}_{E,X} &= \varprojlim (\mathcal{O}_E/\omega^n\mathcal{O}_E)_X && \text{for } E/\mathbb{Q}_{\ell} \text{ finite};\\
\mathcal{O}_{E,X} &= \underset{\substack{ F / \mathbb{Q}_{\ell} \text{ finite,} \\ F \subset  E}}{\colim} \mathcal{O}_{F,X} && \text{for } E/\mathbb{Q}_{\ell} \text{ algebraic};\\
E_X & = \mathcal{O}_{E,X}[l^{-1}]  &&\text{for } E/\mathbb{Q}_{\ell} \text{ algebraic}.\\
\end{align*}
\end{Example}
\begin{Lemma}\label{Lemma: pullback is exact in pro-etale topology}
For any morphism of schemes $X \rightarrow Y$, the canonical map $f^{\#}: f^{-1}\Lambda_Y \rightarrow \Lambda_X$ is an isomorphism. In particular $f^* = f^{-1}$ is exact.
\end{Lemma}
\begin{proof}
This result can be found in \cite[Proposition 4.6]{chough2020pro}. The proof uses the explicit description of $\Lambda_X$ in \cref{Example: explicit description of sheaf of ring in pro-etale topology}. Since pullback commutes with colimits and localization, the main step is to show that $f^{\#}$ is an isomorphism when $\Lambda = \mathcal{O}_E$ for a finite extension $E/\mathbb{Q}_{\ell}$. In this case, $\Lambda_Y =  \varprojlim (\mathcal{O}_E/\omega^n\mathcal{O}_E)_Y $ is represented by the scheme $ \varprojlim (\bigsqcup_{\mathcal{O}_E/\omega^n\mathcal{O}_E} Y)$. Hence, $f^*\Lambda_Y$ is represented by the scheme
\begin{equation*}
X \times_Y \varprojlim (\bigsqcup_{\mathcal{O}_E/\omega^n\mathcal{O}_E} Y) = \varprojlim ( \bigsqcup_{\mathcal{O}_E/\omega^n\mathcal{O}_E} X \times_Y Y) = \varprojlim ( \bigsqcup_{\mathcal{O}_E/\omega^n\mathcal{O}_E} X)
\end{equation*}
and thus $f^*\Lambda_Y = \Lambda_X$.
\end{proof}
Recall that a scheme $X$ is called $\ell$-coprime if $\ell$ is invertible in $\mathcal{O}_X$. Following \cite{bhatt2015pro}, the constructible derived category can be defined as follows.
\begin{Definition}\label{Definition: derived category in proetale site}
Let $X$ be a noetherian $\ell$-coprime scheme.
\begin{enumerate}
\item A sheaf $\mathcal{F} \in \Sh(X_{\proet}, \Lambda_X)$ is called \textit{locally constant of finite presentation} if it is locally isomorphic to $M \otimes_{\Lambda} \Lambda_X$ for a finitely-presented $\Lambda$-module $M$.
\item A sheaf $\mathcal{F} \in \Sh(X_{\proet}, \Lambda_X)$ is called \textit{constructible} if there exists a finite stratification $\{X_i \rightarrow X\}$ such that $\mathcal{F}|_{X_i}$ is locally constant of finite presentation.
\item A complex $\mathcal{F} \in D(X_{\proet}, \Lambda_X)$ is called \textit{constructible} if it is bounded and all its cohomology sheaves are constructible.
\end{enumerate}
We denote by $\Sh_c(X, \Lambda) \subset \Sh(X_{\proet}, \Lambda_X)$ the category of all constructible sheaves. The \textit{constructible derived category} is the full subcategory $D^b_c(X, \Lambda) \subset D(X_{\proet}, \Lambda_X)$ consisting of constructible complexes.
\end{Definition}
\begin{Remark}
The assumptions on $X$ in the definition above can be weakened to topologically noetherian and qcqs. For our purposes, noetherian is good enough since all schemes we will encounter are noetherian.
\end{Remark}
The six functor formalism can also be conveniently described in the pro-étale topology: The functors $\otimes^L, \RIHom$ and $f^*$ preserve constructible complexes and thus they descend to the constructible derived category. The same is true for $f_*$ if $f :X \rightarrow Y$ is of finite type and $Y$ satisfies some mild assumptions (e.g. for $Y$ of finite type over a field or more generally for $Y$ quasi-excellent). If $j : U \hookrightarrow X$ is an open immersion, the pullback functor $j^*$ has a left adjoint $j_!$ which preserves constructible complexes. For a general $f:X \rightarrow Y$ (of finite type with $Y$ quasi-excellent), we define $f_! := \overline{f}_* \circ j_!$ where we factor $f$ as an open immersion followed by a proper map
\begin{equation}\label{eq: factorization of f into open and proper}
X \overset{j}{\hookrightarrow} \overline{X} \overset{\overline{f}}{\rightarrow} Y.
\end{equation}
When restricted to the constructible derived category, the functor $f_!$ admits a right adjoint denoted by $f^!$. These are the usual six functors $f^*,f_*,f_!,f^!, \otimes^L$ and $\RIHom$.
\begin{Remark}
In \cite{bhatt2015pro} a slightly different pullback functor $f_{comp}^*$ is used which is defined as $f^*$ followed by a certain completion operation (\cite[Lemma 6.5.9]{bhatt2015pro}). It turns out that in our situation (i.e. under \cref{Assumption: coefficient ring Lambda}) the two functors agree on the constructible derived category. To see this, recall that $f^*$ is exact by \cref{Lemma: pullback is exact in pro-etale topology}. Using this, it is straightforward to check that $f^*$ preserves constructible complexes. Since constructible complexes are complete (\cite[Proposition 6.8.11(3), Definition 6.5.1]{bhatt2015pro}), it follows that $f^* = f^*_{comp}$ on $D^b_c(X, \Lambda)$.
\end{Remark}
If $\mathbb{F}$ is a finite field of characteristic $\ell$, we denote by $D^b_c(X_{\et}, \mathbb{F})$ the constructible derived category in the étale topology with coefficients in $\mathbb{F}$. The standard compatibility and base change results for the six functors on $D^b_c(X, \Lambda)$ can be deduced from the corresponding results in $D^b_c(X_{\et}, \mathbb{F})$ using the following reduction steps.
\begin{Lemma}\cite{bhatt2015pro}\label{Lemma: reduction to etale}
Let $X$ be noetherian and $\ell$-coprime.
\begin{enumerate}
\item Let $\mathbb{F}$ be a finite field of characteristic $\ell$. Then the canonical morphism of sites $\nu : X_{\proet} \rightarrow X_{\et}$ induces an equivalence of categories
\begin{equation*}
\nu^* : D^b_c(X_{\et}, \mathbb{F} ) \overset{\sim}{\rightarrow}  D^b_c(X, \mathbb{F});
\end{equation*}
\item Let $E/\mathbb{Q}_{\ell}$ be a finite extension and $\kappa$ the residue field of $\mathcal{O}_E$. Then the functor 
\begin{equation*}
\kappa_X \otimes^L_{\mathcal{O}_{E,X}} - : D^b_c(X, \mathcal{O}_E) \rightarrow D^b_c(X,\kappa)
\end{equation*}
is well-defined (i.e. preserves constructibility) and conservative (i.e. reflects isomorphisms);
\item Let $E / \mathbb{Q}_{\ell}$ be algebraic. Then the canonical functor 
\begin{equation*}
\underset{\substack{ F / \mathbb{Q}_{\ell} \text{ finite,} \\ F \subset  E}}{\colim} D^b_c(X, \mathcal{O}_F) \rightarrow D^b_c(X, \mathcal{O}_E)
\end{equation*}
is an equivalence of categories;
\item Let $E / \mathbb{Q}_{\ell}$ be algebraic. Then the canonical functor
\begin{equation*}
D^b_c(X, \mathcal{O}_{E})[\ell^{-1}] \rightarrow  D^b_c(X, E)
\end{equation*}
is an equivalence of categories. Here $D^b_c(X, \mathcal{O}_{E})[\ell^{-1}]$ is the category with the same objects as $D^b_c(X, \mathcal{O}_{E})$ and
\begin{equation*}
\Hom_{D^b_c(X, \mathcal{O}_{E})[\ell^{-1}]} (\mathcal{F}, \mathcal{G}) := \Hom_{D^b_c(X, \mathcal{O}_{E})} (\mathcal{F}, \mathcal{G})[\ell^{-1}].
\end{equation*}
\end{enumerate}
Moreover, the functors in (1)-(4) are compatible with the six functors $f^*,f_*,f_!,f^!,\otimes^L$ and $\RIHom$.
\end{Lemma}

\begin{proof}
This follows from various results in \cite{bhatt2015pro}: (1) follows from \cite[Corollary 5.1.5, 5.1.6]{bhatt2015pro}.  The well-definedness in (2) follows from \cite[Proposition 6.8.11(3), Definition 6.5.1]{bhatt2015pro}. To show that the functor $\kappa_X \otimes^L_{\mathcal{O}_{E,X}} -$ is conservative, it suffices to prove that it reflects zero objects. Let $\mathcal{F} \in D^b_c(X, \mathcal{O}_E)$ with $ \kappa_X \otimes^L_{\mathcal{O}_{E,X}}  \mathcal{F} = 0$. We need to show that $\mathcal{F} = 0$. Let $\mathfrak{m} \subset \mathcal{O}_E$ be the maximal ideal. Then there is a short exact sequence of $\mathcal{O}_E$-modules
\begin{equation*}
0 \rightarrow \mathfrak{m}^n/\mathfrak{m}^{n+1} \rightarrow \mathcal{O}_E/\mathfrak{m}^{n+1} \rightarrow \mathcal{O}_E/ \mathfrak{m}^n \rightarrow 0
\end{equation*}
with $ \mathfrak{m}^n/\mathfrak{m}^{n+1} \cong \kappa$. This induces a distinguished triangle
\begin{equation*}
\kappa_X \otimes^L_{\mathcal{O}_{E, X}} \mathcal{F}  \rightarrow (\mathcal{O}_E/ \mathfrak{m}^{n+1})_X \otimes^L_{\mathcal{O}_{E,X}}  \mathcal{F} \rightarrow (\mathcal{O}_E/\mathfrak{m}^n)_X \otimes^L_{\mathcal{O}_{E,X}}  \mathcal{F} \rightarrow [1].
\end{equation*}
By an induction argument, this implies that $(\mathcal{O}_E/\mathfrak{m}^n)_X \otimes^L_{\mathcal{O}_{E,X}} \mathcal{F}  = 0$ for all $n\ge 1$. Moreover, $\mathcal{F}$ is $\mathfrak{m}$-adically complete by \cite[Proposition 6.8.11(3), Definition 6.5.1]{bhatt2015pro} which means that
\begin{equation*}
\mathcal{F} = \Rlim ((\mathcal{O}_E/\mathfrak{m}^n)_X \otimes^L_{\mathcal{O}_{E,X}} \mathcal{F}) = 0.
\end{equation*}
This completes the proof of (2). The claims in (3) and (4) are \cite[Proposition 6.8.14]{bhatt2015pro}. It remains to prove compatibility with the six functors in (1)-(4). Note for this that the category $\colim_F D^b_c(X, \mathcal{O}_F)$ in (3) inherits the six functors from the $D^b_c(X, \mathcal{O}_F)$ (a standard argument shows that the transition maps in the colimit are compatible with the six functors using that $\mathcal{O}_{F',X} \cong  \mathcal{O}_{F,X}^{\oplus [F':F]} $ as sheaves of $\mathcal{O}_{F,X}$-modules for any finite extension $F'/F$). Similarly, the category $D^b_c(X, \mathcal{O}_{E})[\ell^{-1}]$ in (4) inherits the six functors from $D^b_c(X, \mathcal{O}_{E})$. All functors in (1)-(4) are induced by pullbacks along morphisms of ringed sites (e.g. the functor in (2) is the (derived) pullback along the morphism of ringed sites $(X_{\proet}, \kappa_X) \rightarrow (X_{\proet}, \mathcal{O}_{E,X})$). As such, they commute with $f^*$ and $\otimes^L$ by standard results about ringed sites \cite[\href{https://stacks.math.columbia.edu/tag/0D6D}{Tag 0D6D}, \href{https://stacks.math.columbia.edu/tag/07A4}{Tag 07A4}]{stacks-project}. Since the functors in (1),(3) and (4) are equivalences, they also commute with the corresponding (right) adjoints $f_*$ and $\RIHom$. Similarly, they commute with $j_!$ for $j$ an open immersion which is left adjoint to $j^*$. Hence, they also commute with $f_! = \overline{f}_* \circ j_! $ (see \eqref{eq: factorization of f into open and proper}) and thus also with its right adjoint $f^!$. It remains to show that the functor in (2) commutes with $f_*, \RIHom, f_!$ and $f^!$. This follows from \cite[Lemma 6.5.11(3), 6.7.13, 6.7.14, 6.7.19]{bhatt2015pro}.
\end{proof}
We conclude this section by collecting a few useful properties about pro-étale sheaves. In the pro-étale topology, so-called $w$-contractible affine schemes play a distinguished role.
\begin{Definition}
An affine scheme $U$ is called $w$-contractible if every faithfully-flat weakly-étale map $V\rightarrow U$ has a section.
\end{Definition}
\begin{Lemma}\label{Lemma: good behaviour of w-contractible schemes}
Let $U$ be a $w$-contractible affine scheme.
\begin{enumerate}
\item The global sections functor $\Gamma(U,-)$ is exact on $\Sh(U_{\proet}, \Lambda_U)$;
\item Any locally constant sheaf of finite presentation $\mathcal{F} \in \Sh_c(U, \Lambda)$ is already constant, i.e. $\mathcal{F} \cong M \otimes_{\Lambda} \Lambda_U$ for a finitely-presented $\Lambda$-module $M$;
\item If $R$ is a strictly Henselian local ring (e.g. an algebraically closed field) then $\Spec(R)$ is $w$-contractible.
\end{enumerate}
\end{Lemma}
\begin{proof}
The exactness in (1) is mentioned in the introduction of \cite{bhatt2015pro} (see also \cite[\href{https://stacks.math.columbia.edu/tag/098H}{Tag 098H}, \href{https://stacks.math.columbia.edu/tag/0946}{Tag 0946}]{stacks-project}). Let $\mathcal{F} \in \Sh_c(U, \Lambda)$ be locally constant of finite presentations. Then we can find a weakly-étale cover $V = \bigsqcup V_i \rightarrow U$ with $\mathcal{F}|_V = M \otimes_{\Lambda} \Lambda_V$ for some finitely-presented $\Lambda$-module $M$. Since $U$ is weakly contractible, there is a section $s:U \rightarrow V$ and thus
\begin{equation*}
\mathcal{F} = s^*(\mathcal{F}|_V) = s^*(M \otimes_{\Lambda} \Lambda_V ) = M \otimes_{\Lambda} \Lambda_U.
\end{equation*}
This proves (2). (3) follows from \cite[Theorem 1.8]{bhatt2015pro}.
\end{proof}
\begin{Lemma}\label{Lemma: construcible derived category is closed under direct summands}
The triangulated category $D^b_c(X, \Lambda) \subset D(X_{\proet}, \Lambda_X)$ is closed under direct summands. In particular, $D^b_c(X, \Lambda)$ is idempotent complete.
\end{Lemma}
\begin{proof}
Since taking cohomology commutes with direct sums, it suffices to prove that $\Sh_c(X, \Lambda) \subset \Sh(X_{\proet}, \Lambda_X)$ is closed under direct summands. In fact, by \cite[Lemma 6.8.7, Proposition 6.8.11]{bhatt2015pro} the category $\Sh_c(X, \Lambda)$ is abelian so it is certainly closed under direct summands.
\end{proof}
\begin{Lemma}\label{Lemma: pullback along pro-etale morphism}
If $f:X \rightarrow Y$ is weakly étale, the pullback functor $f^*$ preserves injectives.
\end{Lemma}
\begin{proof}
Compositions of weakly étale maps are weakly étale and morphisms between weakly étale maps are weakly étale. Together with \cref{Lemma: pullback is exact in pro-etale topology} this implies that $(X_{\proet}, \Lambda_X)$ is the localization of the ringed site $(Y_{\proet}, \Lambda_Y)$ at $X \rightarrow Y$. By general results on ringed sites \cite[\href{https://stacks.math.columbia.edu/tag/04IX}{Tag 04IX}]{stacks-project} this implies that $f^*$ has an exact left adjoint $f_!$. In particular, $f^*$ preserves injectives.
\end{proof}
Let $\mathcal{F} \in D^b_c(X, \Lambda)$ and pick a complex of injectives $I_{\mathcal{F}}^{\bullet} \in C^+(\Sh(X_{\proet}, \Lambda_{X}))$ representing $\mathcal{F}$. Let $R_{\mathcal{F}} = \Hom_{dg}^{\bullet}(I_{\mathcal{F}}^{\bullet},I_{\mathcal{F}}^{\bullet})$ be the dg-algebra from \eqref{eq: definition of dg algebra associed to X}. Thanks to the pro-étale formalism, we obtain the following algebraic description of the category $\langle \mathcal{F} \rangle_{\Delta, \inplus}$.
\begin{Corollary}\label{Cor: Equivalence between sheaves and dg-modules for pro-etale}
There is an equivalence of triangulated categories
\begin{equation*}
  \langle  \mathcal{F} \rangle_{D^b_c(X, \Lambda)^{\op},\Delta, \inplus}\cong D_{\perf} (R_{\mathcal{F}}- \dgMod)
\end{equation*}
which send $\mathcal{F}$ to the free dg-module $R_{\mathcal{F}}$.
\end{Corollary}
\begin{proof}
By \cref{Lemma: construcible derived category is closed under direct summands} we have $\langle \mathcal{F} \rangle_{D^b_c(X, \Lambda),\Delta, \inplus} = \langle \mathcal{F} \rangle_{D^+(X_{\proet}, \Lambda_X),\Delta, \inplus}$. The result now follows from \cref{Cor: derived cat generated by X equivalent to perfect dgMod cat}.
\end{proof}
\subsection{The Frobenius action on dg-algebras}\label{Section: The Frobenius action on dg-algebras}
In this section we show that the Frobenius action on $\uHom$ from \eqref{eq: def of uHom on constructible sheaves} is compatible with the dg-techniques from \cref{Section: Derived categories and dg-algebras}. These results are similar to the ones in \cite[Appendix A]{polishchuk2019semiorthogonal} but some of the arguments simplify because the canonical morphism $\Spec(\overline{\mathbb{F}}_q) \rightarrow \Spec(\mathbb{F}_q)$ is a weakly-étale (but not étale).
\medskip

Let $X_0$ be a variety defined over a finite field $\mathbb{F}_q$. Recall that we define $X := X_0 \otimes_{\mathbb{F}_q} \overline{\mathbb{F}}_q$. Let $\pi :X \rightarrow X_0$ be the canonical map and define $\mathcal{F} := \pi^*(\mathcal{F}_0) \in \Sh(X_{\proet}, \Lambda_X)$ for any $\mathcal{F}_0 \in \Sh(X_{0, \proet}, \Lambda_{X_0})$.
\begin{Lemma}\label{Lemma: Pullback to algebraic closure preserves injectives}
The pullback functor $\pi^*:  \Sh(X_{0, \proet}, \Lambda_{X_0}) \rightarrow \Sh(X_{\proet}, \Lambda_X)$ preserves injectives.
\end{Lemma}
\begin{proof}
The morphism $\Spec(\overline{\mathbb{F}}_q) \rightarrow \Spec(\mathbb{F}_q)$ is weakly étale. By base change, $\pi^*: X \rightarrow X_0$ is also weakly étale. Hence, $\pi^*$ preserves injectives by \cref{Lemma: pullback along pro-etale morphism}.
\end{proof}
The geometric $q$-Frobenius $F \in \Gal(\overline{\mathbb{F}}_q / \mathbb{F}_q)$ defines a morphism $\Spec(\overline{\mathbb{F}}_q) \rightarrow \Spec(\overline{\mathbb{F}}_q)$ in the category $\Spec(\mathbb{F}_q)_{\proet}$. Thus, for any $\mathcal{F}_0 \in \Sh( \Spec(\mathbb{F}_q)_{\proet}, \Lambda_{\Spec(\mathbb{F}_q)} )$ there is a corresponding restriction map on sections
\begin{equation*}
F : \Gamma(\Spec(\overline{\mathbb{F}}_q), \mathcal{F}_0 ) \rightarrow \Gamma(\Spec(\overline{\mathbb{F}}_q), \mathcal{F}_0 ). 
\end{equation*}
This map is an isomorphism with inverse given by restriction along the arithmetic Frobenius $F^{-1} \in \Gal(\overline{\mathbb{F}}_q / \mathbb{F}_q)$. Hence, we can equip $\Gamma(\Spec(\overline{\mathbb{F}}_q), \mathcal{F}_0 )$ with the structure of a $\Lambda[F,F^{-1}]$-module. Moreover, for any morphism of sheaves $\mathcal{F}_0 \rightarrow \mathcal{G}_0$, we get a commutative diagram
\begin{equation*}
\begin{tikzcd}
{\Gamma(\Spec(\overline{\mathbb{F}}_q), \mathcal{F}_0 )} \arrow[r, "F"] \arrow[d] & {\Gamma(\Spec(\overline{\mathbb{F}}_q), \mathcal{F}_0 )} \arrow[d] \\
{\Gamma(\Spec(\overline{\mathbb{F}}_q), \mathcal{G}_0 )} \arrow[r, "F"]                                                           & {\Gamma(\Spec(\overline{\mathbb{F}}_q), \mathcal{G}_0 )}.\end{tikzcd}
\end{equation*}
Hence, we obtain a functor
\begin{align*}
\Gamma_{F} : \Sh(\Spec(\mathbb{F}_q)_{\proet}, \Lambda_{\Spec(\mathbb{F}_q)}) & \rightarrow \Lambda[F,F^{-1}]-\modu \\
\mathcal{F}_0 & \mapsto \Gamma(\Spec(\overline{\mathbb{F}}_q), \mathcal{F}_0).
\end{align*}
\begin{Remark}\label{Remark: Gamma Fr extends equivalence between Shc(F_q) and Gal-reps}
Recall from \cref{Section: Frobenius action} that there is an equivalence between the category of constructible sheaves $\Sh_c(\Spec(\mathbb{F}_q), \overline{\mathbb{Q}}_{\ell})$ and the category of finite-dimensional continuous $\Gal(\overline{\mathbb{F}}_q / \mathbb{F}_q)$-representations over $\overline{\mathbb{Q}}_{\ell}$. On $\Sh_c(\Spec(\mathbb{F}_q), \overline{\mathbb{Q}}_{\ell})$ the functor $\Gamma_{F}$ simply corresponds to restricting the $\Gal(\overline{\mathbb{F}}_q / \mathbb{F}_q)$-action to the (geometric) Frobenius element $F \in \Gal(\overline{\mathbb{F}}_q / \mathbb{F}_q)$.
\end{Remark}
Note that
\begin{equation}\label{eq: forgetting Frobenius recovers sections}
\For_{F} \circ \Gamma_{F} (-)= \Gamma(\Spec(\overline{\mathbb{F}}_q), -).
\end{equation}
where
\begin{equation*}
\For_{F}: \Lambda[F,F^{-1}]-\modu \rightarrow \Lambda-\modu
\end{equation*}
is the forgetful functor.
\begin{Lemma}
The functors $\For_{F}, \Gamma_{F}$ and $ \Gamma(\Spec(\overline{\mathbb{F}}_q), -)$ are exact.
\end{Lemma}
\begin{proof}
The forgetful functor $\For_{F}$ is exact and reflects exact sequences. Hence, by \eqref{eq: forgetting Frobenius recovers sections} it suffices to prove that $\Gamma(\Spec(\overline{\mathbb{F}}_q), -)$ is exact. This follows from \cref{Lemma: good behaviour of w-contractible schemes}.
\end{proof}
We define a functor
\begin{equation*}
\uHom(-,-) :  \Sh(X_{0,\proet}, \Lambda_{X_0})^{\op} \times \Sh(X_{0,\proet}, \Lambda_{X_0}) \rightarrow \Lambda[F,F^{-1}]-\modu 
\end{equation*}
via 
\begin{equation*}
(\mathcal{F}_0, \mathcal{G}_0) \mapsto \Gamma_{F}\circ a_* \circ \IHom(\mathcal{F}_0, \mathcal{G}_0).
\end{equation*}
Note that we have
\begin{align*}
\For_{F} (\uHom( \mathcal{F}_0, \mathcal{G}_0)) & = \For_{F} \circ  \Gamma_{F}\circ a_* \circ \IHom( \mathcal{F}_0, \mathcal{G}_0) \\
&\overset{\eqref{eq: forgetting Frobenius recovers sections}}{=}  \Gamma(\Spec(\overline{\mathbb{F}}_q), a_* \IHom( \mathcal{F}_0, \mathcal{G}_0)) \\
&= \Gamma(X, \IHom( \mathcal{F}_0, \mathcal{G}_0)) \\
&= \Hom( \mathcal{F}, \mathcal{G}).
\end{align*}
Hence, $\uHom(\mathcal{F}_0, \mathcal{G}_0)$ is just the vector space $\Hom( \mathcal{F}, \mathcal{G})$ together with a canonical Frobenius action. The composition map
\begin{equation*}
\IHom(\mathcal{G}_0, \mathcal{K}_0) \otimes \IHom(\mathcal{F}_0, \mathcal{G}_0) \rightarrow \IHom(\mathcal{F}_0, \mathcal{K}_0)
\end{equation*}
gives rise to a map
\begin{equation*}
a_*\IHom(\mathcal{G}_0, \mathcal{K}_0) \otimes a_*\IHom(\mathcal{F}_0, \mathcal{G}_0) \rightarrow a_*\IHom(\mathcal{F}_0, \mathcal{K}_0)
\end{equation*}
and thus to a map of $\Lambda[F,F^{-1}]$-modules
\begin{equation}\label{eq: composition in uHom via IHom}
\uHom(\mathcal{G}_0, \mathcal{K}_0) \otimes_{\Lambda} \uHom(\mathcal{F}_0, \mathcal{G}_0) \rightarrow \uHom(\mathcal{F}_0, \mathcal{K}_0).
\end{equation}
Similarly, the unit morphism
\begin{equation*}
\Lambda_{X_0} \rightarrow \IHom(\mathcal{F}_0, \mathcal{F}_0)
\end{equation*}
induces a map
\begin{equation*}
\Lambda_{\Spec(\mathbb{F}_q)} \rightarrow a_*\IHom(\mathcal{F}_0, \mathcal{F}_0)
\end{equation*}
by adjunction and thus a morphism of $\Lambda[F,F^{-1}]$-modules
\begin{equation}\label{eq: unit map in uHom}
\Lambda \rightarrow \uHom(\mathcal{F}_0, \mathcal{F}_0).
\end{equation}
Forgetting the Frobenius action in \eqref{eq: composition in uHom via IHom} and \eqref{eq: unit map in uHom} recovers the standard composition map
\begin{equation*}
\Hom(\mathcal{G}, \mathcal{K}) \otimes_{\Lambda} \Hom(\mathcal{F}, \mathcal{G}) \rightarrow \Hom(\mathcal{F}, \mathcal{K})
\end{equation*}
and unit map
\begin{equation*}
\Lambda \rightarrow \Hom(\mathcal{F}, \mathcal{F}).
\end{equation*}
Similarly, one can define pre- and post-composition maps
\begin{equation}\label{eq: pre and post composition on uHom}
\begin{aligned}
- \circ \alpha : \uHom(\mathcal{F}_0', \mathcal{G}_0) & \rightarrow \uHom(\mathcal{F}_0, \mathcal{G}_0) \\
\beta \circ - : \uHom(\mathcal{F}_0, \mathcal{G}_0) & \rightarrow \uHom(\mathcal{F}_0, \mathcal{G}_0')
\end{aligned}
\end{equation}
for any $\alpha : \mathcal{F}_0 \rightarrow \mathcal{F}_0'$ and $\beta: \mathcal{G}_0 \rightarrow \mathcal{G}_0'$ which recover the maps
\begin{align*}
- \circ \pi^*(\alpha) : & \Hom(\mathcal{F}', \mathcal{G})  \rightarrow \Hom(\mathcal{F}, \mathcal{G}) \\
\pi^*(\beta) \circ - : & \Hom(\mathcal{F}, \mathcal{G})  \rightarrow \Hom(\mathcal{F}, \mathcal{G}') 
\end{align*}
when forgetting the Frobenius action.

We also define a functor on derived categories
\begin{align*}
\RuHom(-,-) :  D(X_{0,\proet}, \Lambda_{X_0})^{\op} \times D(X_{0,\proet}, \Lambda_{X_0}) \rightarrow D(\Lambda[F,F^{-1}]-\modu)
\end{align*}
via
\begin{equation*}
(\mathcal{F}_0, \mathcal{G}_0)  \mapsto  \Gamma_{F}\circ a_* \circ \RIHom(\mathcal{F}_0, \mathcal{G}_0)
\end{equation*}
(recall that we adhere to the derived convention, i.e. $\Gamma_{F}$ and $a_*$ are understood to be derived when applied to complexes). Moreover, we set
\begin{equation*}
\uHom^i(\mathcal{F}_0, \mathcal{G}_0) := H^i(\RuHom(\mathcal{F}_0, \mathcal{G}_0)) = \Gamma_{F}( H^i(a_*\RIHom(\mathcal{F}_0, \mathcal{G}_0))).
\end{equation*}
\begin{Remark}
This definition of $\uHom^i$ extends our previous definition \eqref{eq: def of uHom on constructible sheaves} for constructible complexes to arbitrary complexes in $D(X_{0,\proet}, \Lambda_{X_0})$ (see also \cref{Remark: Gamma Fr extends equivalence between Shc(F_q) and Gal-reps}).
\end{Remark}
For any two chain complexes $\mathcal{F}_0^{\bullet}, \mathcal{G}_0^{\bullet} \in C(\Sh(X_{0, \proet}, \Lambda_{X_0}))$ we define a chain complex
\begin{equation*}
\uHom_{dg}^{\bullet}(\mathcal{F}_0^{\bullet}, \mathcal{G}_0^{\bullet})\in C(\Lambda[F,F^{-1}]-\modu).
\end{equation*}
Explicitly, this complex is given by
\begin{equation*}
\uHom_{dg}^i (\mathcal{F}_0^{\bullet}, \mathcal{G}_0^{\bullet}) := \prod_{l+m=i} \uHom(\mathcal{F}_0^{-l}, \mathcal{G}_0^m)
\end{equation*}
with differential $d_{\uHom^{\bullet}_{dg}(\mathcal{F}_0^{\bullet},\mathcal{G}_0^{\bullet})} (f) = d_{\mathcal{G}_0^{\bullet}} \circ f - (-1)^{|f|} f \circ d_{\mathcal{F}_0^{\bullet}}$ (where we use \eqref{eq: pre and post composition on uHom}). The composition map on $\uHom(-,-)$ from \eqref{eq: composition in uHom via IHom} can be applied component-wise to define a composition map
\begin{equation*}
\uHom_{dg}^i(\mathcal{G}_0^{\bullet},\mathcal{K}_0^{\bullet}) \otimes_{\Lambda}   \uHom_{dg}^j(\mathcal{F}_0^{\bullet},\mathcal{G}_0^{\bullet})\rightarrow \uHom_{dg}^{i+j}(\mathcal{F}_0^{\bullet},\mathcal{K}_0^{\bullet}).
\end{equation*}
This equips $ \uHom^{\bullet}_{dg}(\mathcal{F}_0^{\bullet},\mathcal{F}_0^{\bullet})$ with the structure of a dg-algebra over $\Lambda$ together with a dg-algebra automorphism induced by the Frobenius action. Forgetting this Frobenius action recovers the dg-algebra $\Hom^{\bullet}_{dg}(\mathcal{F}^{\bullet},\mathcal{F}^{\bullet})$.
\begin{Lemma}\label{Lemma: dg uHom computes derived uHom}
Let $\mathcal{F}_0, \mathcal{G}_0 \in D^+(X_{0, \proet}, \Lambda_{X_0})$ and pick bounded below complexes of injectives $I_{\mathcal{F}_0}^{\bullet}, I_{\mathcal{G}_0}^{\bullet} \in C^+(\Sh(X_{0, \proet}, \Lambda_{X_0}))$ representing $\mathcal{F}_0$ and $\mathcal{G}_0$. Then
\begin{equation*}
H^i( \uHom^{\bullet}_{dg} (I_{\mathcal{F}_0}^{\bullet}, I_{\mathcal{G}_0}^{\bullet} )) = \uHom^i(\mathcal{F}_0, \mathcal{G}_0)
\end{equation*}
as $\Lambda[F,F^{-1}]$-modules.
\end{Lemma}
\begin{proof}
By definition, $\RIHom(\mathcal{F}_0, \mathcal{G}_0)$ is represented by the complex $\IHom_{dg}^{\bullet} (I_{\mathcal{F}_0}^{\bullet},I_{\mathcal{G}_0}^{\bullet} )$ defined as
\begin{equation*}
\IHom_{dg}^i (I_{\mathcal{F}_0}^{\bullet},I_{\mathcal{G}_0}^{\bullet} ) := \prod_{l+m=i} \IHom(I_{\mathcal{F}_0}^{-l}, I_{\mathcal{G}_0}^m)
\end{equation*}
with differential $d_{\IHom^{\bullet}_{dg}(I_{\mathcal{F}_0}^{\bullet},I_{\mathcal{G}_0}^{\bullet})} (f) = d_{I_{\mathcal{G}_0}^{\bullet}} \circ f - (-1)^{|f|} f \circ d_{I_{\mathcal{F}_0}^{\bullet}}$. Note that $\IHom(I_{\mathcal{F}_0}^{-l}, I_{\mathcal{G}_0}^m)$  is acyclic for $a_*$ for all $l,m \in \mathbb{Z}$ (see \cite[V-(4.10), V-(5.2)]{SGA4}). Hence, $\IHom_{dg}^i (\mathcal{F}_0^{\bullet}, \mathcal{G}_0^{\bullet})$ is acyclic for $a_*$ for all $i \in \mathbb{Z}$. Since $\Gamma_{F}$ is exact, $\IHom_{dg}^i (\mathcal{F}_0^{\bullet}, \mathcal{G}_0^{\bullet})$ is also acyclic for $\Gamma_{F} \circ a_*$ for all $i \in \mathbb{Z}$. In particular, $\RuHom(\mathcal{F}_0, \mathcal{G}_0) = \Gamma_{F}\circ a_* \circ \RIHom(\mathcal{F}_0, \mathcal{G}_0)$ is represented by the complex obtained by applying $\Gamma_{F} \circ a_*$ to each component in $\IHom_{dg}^{\bullet} (I_{\mathcal{F}_0}^{\bullet},I_{\mathcal{G}_0}^{\bullet} )$. This is precisely the complex $\uHom_{dg}^{\bullet} (I_{\mathcal{F}_0}^{\bullet},I_{\mathcal{G}_0}^{\bullet} )$. Hence,
\begin{equation*}
H^i(\uHom_{dg}^{\bullet} (I_{\mathcal{F}_0}^{\bullet},I_{\mathcal{G}_0}^{\bullet} )) = H^i(\RuHom(\mathcal{F}_0, \mathcal{G}_0)) = \uHom^i(\mathcal{F}_0, \mathcal{G}_0).
\end{equation*}
\end{proof}
The constructions of this section can be summarized as follows.
\begin{Proposition}\label{Prop: Frobenius action can be lifted to the dg-algebra for an injective resolution}
For any $\mathcal{F}_0 \in D^b_c(X_0, \Lambda)$, there is a bounded below complex of injectives $I_{\mathcal{F}}^{\bullet} \in C^+(\Sh(X_{\proet}, \Lambda_X))$ representing $\mathcal{F}$ and dg-algebra automorphism $F$ of $R_{\mathcal{F}}= \Hom_{dg}^{\bullet}(I_{\mathcal{F}}^{\bullet}, I_{\mathcal{F}}^{\bullet})$ such that the action of $H^i(F)$ on $H^i(R_{\mathcal{F}}) = \Hom^i(\mathcal{F}, \mathcal{F})$ is the Frobenius action coming from $\uHom^i(\mathcal{F}_0, \mathcal{F}_0)$.
\end{Proposition}
\begin{proof}
Pick a complex of injectives $I_{\mathcal{F}_0}^{\bullet} \in C^+(\Sh(X_{0, \proet}, \Lambda_{X_0}))$ representing $\mathcal{F}_0$. Then $I_{\mathcal{F}}^{\bullet} := \pi^*(I_{\mathcal{F}_0}^{\bullet})$ is a complex of injectives representing $\mathcal{F}$ by \cref{Lemma: Pullback to algebraic closure preserves injectives}. The dg-algebra $\uHom_{dg}^{\bullet}(I_{\mathcal{F}_0}^{\bullet}, I_{\mathcal{F}_0}^{\bullet})$ is just $R_{\mathcal{F}} =  \Hom_{dg}^{\bullet}(I_{\mathcal{F}}^{\bullet}, I_{\mathcal{F}}^{\bullet})$ equipped with a dg-algebra automorphism $F$. By \cref{Lemma: dg uHom computes derived uHom}, we have $H^i(\uHom_{dg}^{\bullet}(I_{\mathcal{F}_0}^{\bullet}, I_{\mathcal{F}_0}^{\bullet})) = \uHom^i(\mathcal{F}_0, \mathcal{F}_0)$. Hence, $H^i(F)$ is the canonical Frobenius action on $\uHom^i(\mathcal{F}_0, \mathcal{F}_0)$.
\end{proof}
\subsection{Formality for the Springer category}
Let $\mu_i: \tilde{V}^i \rightarrow V$ be a finite collection of morphisms of Springer type defined over $\overline{\mathbb{F}}_q$ and let $\Sbf = \bigoplus_{i \in I} \Sbf^i \in D^b_c(V, \overline{\mathbb{Q}}_{\ell})$ be the associated Springer sheaf. Recall from \eqref{eq: definition of springer category} that the Springer category is defined as
\begin{equation*}
D_{\Spr}(V, \overline{\mathbb{Q}}_{\ell}):= \langle X_1, ...,X_n \rangle_{\Delta} \subset D^b_c(V, \overline{\mathbb{Q}}_{\ell})
\end{equation*} 
where the $X_i$ are the simple perverse sheaves appearing in $\Sbf$.
 \begin{Lemma}\label{Cor: Springer category is generated by S under direct summands}
We have $D_{\Spr}(V, \overline{\mathbb{Q}}_{\ell}) = \langle \Sbf \rangle_{ \Delta, \inplus}$.
\end{Lemma}
\begin{proof}
Since $X_1, ..., X_n \in \langle \Sbf \rangle_{\Delta, \inplus}$, we have $D_{\Spr}(V, \overline{\mathbb{Q}}_{\ell}) \subset \langle \Sbf \rangle_{ \Delta, \inplus}$. Hence, it suffices to show that $D_{\Spr}(V, \overline{\mathbb{Q}}_{\ell})$ is closed under direct summands. Let $\Perv_{\Spr}(V)  \subset \Perv(V)$ be the Serre subcategory generated by the $X_1,...,X_n$. Then by a standard dévissage argument, the Springer category can be described as
\begin{equation*}
D_{\Spr}(V, \overline{\mathbb{Q}}_{\ell}) = \{ \mathcal{F} \in D^b_c(V, \overline{\mathbb{Q}}_{\ell}) \mid {}^pH^i(\mathcal{F}) \in \Perv_{\Spr}(V) \text{ for all } i \in \mathbb{Z}\}.
\end{equation*}
Note that as a Serre subcategory, $\Perv_{\Spr}(V)$ is closed under direct summands in $\Perv(V)$. Hence, $D_{\Spr}(V, \overline{\mathbb{Q}}_{\ell})$ is also closed under direct summands.
\end{proof} 
We obtain the following general formality result for the Springer category.
\begin{Theorem}\label{Theorem: Springer formality}
There is an equivalence of triangulated categories
\begin{equation*}
D_{\Spr}(V, \overline{\mathbb{Q}}_{\ell})^{\op} \cong D_{\perf}(\Hom^*(\Sbf,\Sbf)-\dgMod)
\end{equation*}
which identifies $\Sbf$ with the free dg-module $\Hom^*(\Sbf, \Sbf)$.
\end{Theorem}
\begin{proof}
Combining \cref{Cor: Equivalence between sheaves and dg-modules for pro-etale} and \cref{Cor: Springer category is generated by S under direct summands} we get an equivalence
\begin{equation*}
D_{\Spr}(V, \overline{\mathbb{Q}}_{\ell})^{\op} = \langle \Sbf \rangle_{ \Delta, \inplus} \cong  D_{\perf}(R_{\Sbf}-\dgMod)
\end{equation*}
which send $\Sbf$ to the free dg-module $ R_{\Sbf}$. By \cref{Cor: Purity of the Ext algebra} and \cref{Prop: Frobenius action can be lifted to the dg-algebra for an injective resolution} the dg-algebra $R_{\Sbf}$ can be equipped with a dg-algebra automorphism $F: R_{\Sbf} \rightarrow R_{\Sbf}$ such that $H^i(F)$ acts by multiplication with $q^{\tfrac{i}{2}}$. This implies that $R_{\Sbf}$ is formal by \cref{Lemma: purity implies formality}. In particular, by \eqref{eq: equivalence of categories when Rx is formal} we get an equivalence
\begin{equation*}
D_{\Spr}(V, \overline{\mathbb{Q}}_{\ell})^{\op} \cong D_{\perf}(R_{\Sbf}-\dgMod) \cong D_{\perf}(\Hom^*(\Sbf,\Sbf)-\dgMod)
\end{equation*}
which sends $\Sbf$ to the free dg-module $\Hom^*(\Sbf,\Sbf)$.
\end{proof}

\section{De $\mathbb{F}$ \`{a} $\mathbb{C}$}\label{Section: De F a C}
In this section we prove a formality result for Springer sheaves on varieties over a field of characteristic $0$ by reduction to the positive characteristic case. This is a standard application of the "De $\mathbb{F}$ a $\mathbb{C}$" technique from \cite[§6]{beilinson2018faisceaux}. We will explain how these arguments work in the pro-étale setting. If $X_A$ is a scheme defined over a ring $A$ and $A \rightarrow R$ is a ring homomorphism, we denote by $X_R$ the base change of $X_A$ to $R$. Similarly, if $S \rightarrow \Spec(A)$ is a morphism of schemes, we denote by $X_S$ the corresponding base change to $S$. Furthermore, for $\mathcal{F} \in D^b_c(X_A, \overline{\mathbb{Q}}_{\ell})$, we denote by $\mathcal{F}_R \in D^b_c(X_R, \overline{\mathbb{Q}}_{\ell})$ (resp. $\mathcal{F}_S \in D^b_c(X_S, \overline{\mathbb{Q}}_{\ell})$) the corresponding pullback to $X_R$ (resp. $X_S$).

\subsection{Generic base change}
The main tool that we need to compare constructible sheaves on varieties in characteristic $p$ and characteristic $0$ is the generic base change theorem. Let us first explain what we mean by generic base change.
\begin{Definition}
Let $S$ be a scheme and $f:X\rightarrow Y$ a morphism of $S$-schemes. Let $\mathcal{F}, \mathcal{G} \in D(X_{\et}, \Lambda)$ for a noetherian ring $\Lambda$ (resp. $\mathcal{F}, \mathcal{G} \in D(X_{\proet}, \Lambda_X)$ for $\Lambda $ as in \cref{Assumption: coefficient ring Lambda}).
\begin{enumerate}
\item We say that the formation $f_*\mathcal{F}$ commutes with generic base change if there is a dense open subscheme $U \subset S$ such that for each morphism of schemes $g:S' \rightarrow U \subset S$ with corresponding pullback diagram
\begin{equation*}
\begin{tikzcd}
X_{S'} \arrow[d, "f'"] \arrow[r, "g''"] & X \arrow[d, "f"] \\
Y_{S'} \arrow[d] \arrow[r, "g'"]        & Y \arrow[d]      \\
S' \arrow[r, "g"]                           & S     
\end{tikzcd}
\end{equation*}
the canonical map $(g')^*f_* \mathcal{F} \rightarrow f'_*(g'')^*\mathcal{F}$ is an isomorphism.
\item We say that the formation $\RIHom(\mathcal{F}, \mathcal{G})$ commutes with generic base change if there is a dense open subscheme $U \subset X$ such that for each morphism of schemes $g: S' \rightarrow U \subset S$ the canonical map $g^*\RIHom(\mathcal{F}, \mathcal{G}) \rightarrow \RIHom(g^*\mathcal{F}, g^*\mathcal{G})$ is an isomorphism.
\end{enumerate}
\end{Definition}
There is the following generic base change theorem for étale sheaves.
\begin{Theorem}\label{Theorem: generic base change on etale site}\cite{SGA4-1/2}
Let $S$ be a noetherian scheme, $f:X \rightarrow Y$ a morphism of $S$-schemes of finite type and $\Lambda$ a noetherian ring annihilated by an integer invertible in $\mathcal{O}_S$. Then for any $\mathcal{F}, \mathcal{G} \in D^b_c(X_{\et}, \Lambda)$ the formations $f_* \mathcal{F}$ and $\RIHom(\mathcal{F}, \mathcal{G})$ commute with generic base change. 
\end{Theorem}
\begin{proof}
This is proved for constructible sheaves $\mathcal{F} \in \Sh_c(X_{\et}, \Lambda)$ in \cite[Th. finitude 1.9, 2.10]{SGA4-1/2}. The case of an arbitrary complex $\mathcal{F} \in D^b_c(X_{\et}, \Lambda)$ can be deduced from this by a dévissage argument.
\end{proof}
We get a similar result for non-torsion coefficients in the pro-étale topology. For the rest of this section, we assume that $\Lambda$ is as in \cref{Assumption: coefficient ring Lambda}.
\begin{Lemma}\label{Lemma: generic base change for proetale site}
Let $S$ be a noetherian, quasi-excellent, $\ell$-coprime scheme and $f:X \rightarrow Y$ a morphism of $S$-schemes of finite type. Then for any $\mathcal{F} ,\mathcal{G} \in D^b_c(X_{\proet}, \Lambda_X)$ the formations $f_* \mathcal{F}$ and $\RIHom(\mathcal{F}, \mathcal{G})$ commute with generic base change.
\end{Lemma}
\begin{proof}
Using the reduction steps for constructible sheaves from \cref{Lemma: reduction to etale} this can be reduced to the étale generic base change theorem (\cref{Theorem: generic base change on etale site}). The assumption that $S$ is quasi-excellent ensures that $f_* \mathcal{F}$ is constructible.
\end{proof}
The following lemma can also be deduced from \cite[Lemme 6.1.9]{beilinson2018faisceaux} or \cite[A.5]{polishchuk2019semiorthogonal}.
\begin{Lemma}\label{Lemma: lifting triangulated category generated by F}
Let $X_A$ be a scheme of finite type over a finitely-generated $\mathbb{Z}[\ell^{-1}]$-algebra $A \subset \overline{\mathbb{Q}}_{\ell}$ and $s \in \Spec(A)$ a closed point. For any $\mathcal{F}_A \in D^b_c(X_A, \Lambda)$ there exists a strictly Henselian local ring $R$ such that
\begin{enumerate}[label=(\roman*)]
\item $A \subset R \subset \overline{\mathbb{Q}}_{\ell}$;
\item the residue field of $R$ is the algebraic closure $\overline{s}$ of $s$;
\item Pulling back along the morphisms in the diagram $X_{\overline{\mathbb{Q}}_{\ell}} \overset{u}{\leftarrow} X_R \overset{v}{\rightarrow} X_s$ induces equivalences of categories
\begin{equation*}
\langle \mathcal{F}_{\overline{\mathbb{Q}}_{\ell}} \rangle_{D^b_c(X_{\overline{\mathbb{Q}}_{\ell}}, \Lambda), \Delta, \inplus} \cong \langle \mathcal{F}_R \rangle_{D^b_c(X_R, \Lambda), \Delta, \inplus}  \cong \langle \mathcal{F}_{\overline{s}} \rangle_{D^b_c(X_{\overline{s}}, \Lambda), \Delta, \inplus}.
\end{equation*}
\end{enumerate}
\end{Lemma}
\begin{proof}
By constructibility, we may enlarge $A$ so that the cohomology sheaves of $a_*\RIHom(\mathcal{F}_A, \mathcal{F}_A)$ are locally constant of finite presentation, i.e. each $H^i(a_*\RIHom(\mathcal{F}_A, \mathcal{F}_A))$ is locally isomorphic to $M_i \otimes_{\Lambda} \Lambda_{\Spec(A)}$ for some finitely-presented $\Lambda$-module $M_i$. By \cref{Lemma: generic base change for proetale site} we may enlarge $A$ further so that $a_*\RIHom(\mathcal{F}_A, \mathcal{G}_A)$ commutes with arbitrary base change, i.e. for any $g: S \rightarrow \Spec(A)$, we have
\begin{equation*}
g^*a_*\RIHom(\mathcal{F}_A, \mathcal{F}_A) = a_*\RIHom(\mathcal{F}_S, \mathcal{F}_S).
\end{equation*}
Let $S$ be a $w$-contractible affine scheme. Then any locally constant sheaf of finite presentation on $S_{\proet}$ is constant by \cref{Lemma: good behaviour of w-contractible schemes} and hence
\begin{equation}\label{eq: cohomology of RuHom is constand}
\begin{aligned}
H^i(a_*\RIHom(\mathcal{F}_S, \mathcal{F}_S)) & = H^i(g^*a_*\RIHom(\mathcal{F}_A, \mathcal{F}_A))\\
&  =g^* H^i(a_*\RIHom(\mathcal{F}_A, \mathcal{F}_A)) \\
& =  M_i \otimes_{\Lambda} \Lambda_S.
\end{aligned}
\end{equation}
Note that the global sections functor $\Gamma_S = \Gamma(S,-)$ is exact by \cref{Lemma: good behaviour of w-contractible schemes}. Hence, for $S$ affine, connected and $w$-contractible, we get
\begin{equation}\label{eq: hom space is just Mi}
\begin{aligned}
\Hom^i(\mathcal{F}_S, \mathcal{F}_S) &= H^i(\Gamma_S \circ a_* \circ \RIHom(\mathcal{F}_S, \mathcal{F}_S)) \\
&= \Gamma_S \circ H^i( a_*\RIHom(\mathcal{F}_S, \mathcal{F}_S))  \\
& = \Gamma_S (M_i \otimes_{\Lambda} \Lambda_S) \\
&= M_i.
\end{aligned}
\end{equation}
Here the first equality follows from general results about sites (see \cite[V-(4.10), V-(5.2)]{SGA4}). The second uses exactness of $\Gamma_S$ and the third follows from \eqref{eq: cohomology of RuHom is constand}. To see the last equality, pick a presentation $\Lambda^n \rightarrow \Lambda^m \rightarrow M_i \rightarrow 0$. This induces a presentation $ \Lambda_S^n \rightarrow  \Lambda_S^m \rightarrow M_i \otimes_{\Lambda}  \Lambda_S \rightarrow 0$. Note that $\Gamma_S(\Lambda_S) = \Lambda$ by \eqref{eq: constant sheaf evaluated at quasi compact covers} and the fact that $S$ is connected. Since $\Gamma_S$ is exact, this implies $\Gamma_S (M_i \otimes_{\Lambda} \Lambda_S) = \coker( \Lambda^n \rightarrow \Lambda^m) =  M_i$. Now let $R$ be any strictly Henselian local ring satisfying (i) and (ii) (it is explained in \cite[p.156]{beilinson2018faisceaux} that such a ring always exists and one can even assume that $R$ is also a discrete valuation ring). Then the schemes $\Spec(\overline{\mathbb{Q}}_{\ell}), \Spec(R), $ and $\overline{s}$ are all connected $w$-contractible affine schemes by \cref{Lemma: good behaviour of w-contractible schemes}. Using \eqref{eq: hom space is just Mi}, we get
\begin{equation*}
\Hom^i( \mathcal{F}_{\overline{\mathbb{Q}}_{\ell}}, \mathcal{F}_{\overline{\mathbb{Q}}_{\ell}})  = \Hom^i( \mathcal{F}_R, \mathcal{F}_R) = \Hom^i( \mathcal{F}_{\overline{s}}, \mathcal{F}_{\overline{s}}) = M_i
\end{equation*}
for all $i \in \mathbb{Z}$. By \cref{Lemma: idempotent complete Beilinson's Lemma} this implies that there are equivalences of triangulated categories 
\begin{equation*}
\langle \mathcal{F}_{\overline{\mathbb{Q}}_{\ell}} \rangle_{D^b_c(X_{\overline{\mathbb{Q}}_{\ell}}, \Lambda), \Delta, \inplus} \cong \langle \mathcal{F}_R \rangle_{D^b_c(X_R, \Lambda), \Delta, \inplus}  \cong \langle \mathcal{F}_{\overline{s}} \rangle_{D^b_c(X_{\overline{s}}, \Lambda), \Delta, \inplus}.
\end{equation*}
\end{proof}

\subsection{Lifting Springer sheaves}
Any variety $X$ over $\overline{\mathbb{Q}}_{\ell}$ can be defined over a finitely-generated $\mathbb{Z}[\ell^{-1}]$-algebra $A \subset \overline{\mathbb{Q}}_{\ell}$, i.e. there exists a scheme $X_A$ (of finite type) over $\Spec(A)$ such that $X = X_{A} \otimes_{A} \overline{\mathbb{Q}}_{\ell}$. Choosing a maximal ideal of $A$ gives rise to a ring homomorphism $A \rightarrow \mathbb{F}_q$ for a finite field $\mathbb{F}_q$. Hence, it makes sense to consider $X_{\mathbb{F}_q}$ relating $X$ to a scheme over a field of positive characteristic. We first check that this procedure preserves morphisms of Springer type.
\medskip

Let $G$ be a connected reductive group over $\overline{\mathbb{Q}}_{\ell}$. There always exists a split reductive group scheme $G_{\mathbb{Z}}$ over $\mathbb{Z}$ with $G = G_{\mathbb{Z}} \otimes_{\mathbb{Z}} \overline{\mathbb{Q}}_{\ell}$. Fix a Borel subgroup $B_{\mathbb{Z}} \subset G_{\mathbb{Z}}$ and let $B := B_{\mathbb{Z}} \otimes_{\mathbb{Z}} \overline{\mathbb{Q}}_{\ell} $ be the corresponding Borel subgroup of $G$. Let $V$ be a $G$-representation and $\{V^i \subset V \mid i \in I\}$ a finite collection of $B$-stable subspaces. Recall that a representation of an affine group scheme $H$ is the same as a comodule over the Hopf algebra $\mathcal{O}(H)$.
\begin{Lemma}\label{Lemma: reduction mod p preserves morphisms of Springer type}
There is a finitely generated $\mathbb{Z}[\ell^{-1}]$-algebra $A \subset \overline{\mathbb{Q}}_{\ell}$, a $G_A$-representation $V_A$ and $B_A$-stable $A$-submodules $V_A^i \subset V_A$ which recover the $V^i \subset V$ when extending scalars to $\overline{\mathbb{Q}}_{\ell}$.
\end{Lemma}
\begin{proof}
Pick a basis $v_1, ..., v_n \in V$. Then we can write the comultiplication on the $v_i$ as
\begin{equation*}
\Delta(v_i) = \sum_{j} v_j \otimes f_{ij}
\end{equation*}
for some $f_{ij} \in \mathcal{O}(G)$. Since $\mathcal{O}(G) = \mathcal{O}(G_{\mathbb{Z}}) \otimes_{\mathbb{Z}} \overline{\mathbb{Q}}_{\ell}$, we can write each $f_{ij}$ as 
\begin{equation*}
f_{ij} = \sum_k g_{ijk} \otimes c_{ijk}
\end{equation*}
for some $g_{ijk} \in \mathcal{O}(G_{\mathbb{Z}})$ and $c_{ijk} \in \overline{\mathbb{Q}}_{\ell}$. Let $A \subset \overline{\mathbb{Q}}_{\ell}$ be the $\mathbb{Z}[\ell^{-1}]$-algebra generated by all the $c_{ijk}$. Then $V_A := \Span_A \{v_1, ..., v_n\}$ is an $\mathcal{O}(G_A)$-comodule and thus a $G_A$-representation which recovers $V$ when extending scalars to $\overline{\mathbb{Q}}_{\ell}$. We can use the same argument (potentially enlarging $A$ further) to show that there are $B_A$-representations $V^i_A$ which are spanned by finitely many vectors $w^i_{j} \in V^i$ and recover $V^i$ when extending scalars to $\overline{\mathbb{Q}}_l$. After enlarging $A$ further, we may assume that $V_A$ contains all the $w^i_j$ and thus $V^i_A \subset V_A$.
\end{proof}
Consider the variety $\tilde{V}^i_A := G_A \times^{B_A} V^i_A$ and the canonical map $\mu_i: \tilde{V}^i_A \rightarrow V_A$ induced by the action map $G_A \times V^i_A \rightarrow V_A$. For any ring homomorphism $A \rightarrow R$, we get a pullback diagram
\begin{equation*}
\begin{tikzcd}
 \tilde{V}^i_R \arrow[r, "u"]\arrow[d,"\mu_i"]  &  \tilde{V}^i_A \arrow[d,"\mu_i"]    \\
 V_R \arrow[r,"u"]& V_A           
\end{tikzcd}
\end{equation*}
where $\tilde{V}^i_R = G_R \times^{B_R} V^i_R $. We define
\begin{align*}
\Sbf^i_R & := (\mu_i )_! \textbf{1}_{\tilde{V}^i_R } [\dim \tilde{V}^i] \in D^b_c(V_R, \overline{\mathbb{Q}}_{\ell})  \\
\Sbf_R & := \bigoplus_{i \in I} \Sbf^i_R \in D^b_c(V_R, \overline{\mathbb{Q}}_{\ell}).
\end{align*}
By proper base change we have $u^*(\Sbf^i_A) = \Sbf^i_R$ and $u^*(\Sbf_A) = \Sbf_R$. Moreover, if $R = \overline{k}$ is an algebraically closed field, the morphism $\mu_i : \tilde{V}^i_{\overline{k}} \rightarrow V_{\overline{k}}$ is of Springer type and $\Sbf^i_{\overline{k}}$ is the associated Springer sheaf (c.f. \cref{Section: morphisms of Springer type}). Here we use that
\begin{equation*}
\dim \tilde{V}^i = \dim G + \dim V^i - \dim B =  \dim G_{\overline{k}} + \dim V^i_{\overline{k}} - \dim B_{\overline{k}} = \dim \tilde{V}^i_{\overline{k}}.
\end{equation*}
\begin{Corollary}\label{Corollary: Springer formality over C}
Let $\mu_i : \tilde{V}^i \rightarrow V$ be a finite collection of morphisms of Springer type over $\overline{\mathbb{Q}}_{\ell}$. Then there is an equivalence of triangulated categories
\begin{equation*}
D_{\Spr}(V, \overline{\mathbb{Q}}_{\ell})^{\op} \cong  D_{\perf}( \Hom^*(\Sbf, \Sbf)-\dgMod) 
\end{equation*}
which identifies $\Sbf$ with the free dg-module $\Hom^*(\Sbf, \Sbf)$.
\end{Corollary}
\begin{proof}
Let $A$, $\tilde{V}^i_A, V_A,...$ be as in \cref{Lemma: reduction mod p preserves morphisms of Springer type} and pick a closed point $s \rightarrow \Spec(A)$ with algebraic closure $\overline{s}$. Note that since $A$ is finitely-generated over $\mathbb{Z}[\ell^{-1}]$, we have $s = \Spec(\mathbb{F}_q)$ for a finite field $\mathbb{F}_q$ with $\ell$ invertible in $\mathbb{F}_q$. By \cref{Lemma: lifting triangulated category generated by F} there is an equivalence of triangulated categories $D_{\Spr}(V, \overline{\mathbb{Q}}_{\ell}) \cong D_{\Spr}(V_{\overline{s}}, , \overline{\mathbb{Q}}_{\ell})$ which identifies $\Sbf$ with $\Sbf_{\overline{s}}$. In particular, we get $\Hom^*(\Sbf, \Sbf) \cong \Hom^*(\Sbf_{\overline{s}}, \Sbf_{\overline{s}}) $. Hence, by \cref{Theorem: Springer formality}, we have
\begin{equation*}
D_{\Spr}(V, \overline{\mathbb{Q}}_{\ell})^{\op} \cong D_{\Spr}(V_{\overline{s}}, \overline{\mathbb{Q}}_{\ell})^{\op} \cong D_{\perf}( \Hom^*(\Sbf_{\overline{s}}, \Sbf_{\overline{s}})- \dgMod) \cong D_{\perf}(\Hom^*(\Sbf,  \Sbf)-\dgMod).
\end{equation*}
\end{proof}

\section{A derived Deligne-Langlands correspondence}\label{Section: Derived Deligne-Langlands}
Let $G$ be a connected reductive group over $\overline{\mathbb{Q}}_{\ell}$($\cong \mathbb{C}$) with simply connected derived subgroup and let $(X^*,\Phi, X_*,\Phi^{\vee})$ be the associated root datum. Fix a torus and a Borel subgroup $T\subset B \subset G$. Let $\Pi \subset \Phi$ be the associated set of simple roots and $W$ the Weyl group. We recall the definition of the affine Hecke algebra in its Bernstein presentation.
\begin{Definition}
The affine Hecke algebra $\mathcal{H}^{\aff}$ of $G$ is the $\overline{\mathbb{Q}}_{\ell}[q,q^{-1}]$-algebra with generators $\{T_w, \theta_x \mid   w \in W,  x \in X^*\}$ and relations
\begin{itemize}
\item[] \makebox[7cm]{$(T_{s_{\alpha}}+1)(T_{s_{\alpha}}-q) = 0 $ \hfill}  $\alpha \in \Pi $;
\item[] \makebox[7cm]{$ T_wT_{w'} = T_{ww'} $ \hfill} $w,w'\in W$ with $l(w)+l(w') = l(ww')$;
\item[] \makebox[7cm]{$ \theta_x\theta_{x'} = \theta_{x+x'}$ \hfill}  $x,x' \in X^*$;
\item[] \makebox[7cm]{$ T_{s_{\alpha}} \theta_{s_{\alpha}(x)} - \theta_x T_{s_{\alpha}} = (1-q) \frac{\theta_x - \theta_{s_{\alpha}(x)}}{ 1 - \theta_{-\alpha}} $ \hfill}  $\alpha \in \Pi$.
\end{itemize}
\end{Definition}
We collect a few well-known algebraic properties of the affine Hecke algebra which can be found in \cite{lusztig1989affine,chriss2009representation}.
\begin{Lemma}\label{Lemma: algebraic properties of affine Hecke algebra}
\begin{enumerate}
\item $\mathcal{H}^{\aff}$ is a free $\overline{\mathbb{Q}}_{\ell}[q,q^{-1}]$-module with basis $\{\theta_xT_w \mid x\in X^*, w \in W\}$;
\item The $\theta_x$ ($x \in X^*$) span a subalgebra of $\mathcal{H}^{\aff}$ isomorphic to the group algebra $\overline{\mathbb{Q}}_{\ell}[q,q^{-1}][X^*]$;
\item The center of the affine Hecke algebra is $ Z(\mathcal{H}^{\aff}) = \overline{\mathbb{Q}}_{\ell}[q,q^{-1}][X^*]^W$;
\item $\mathcal{H}^{\aff}$ is a free $Z(\mathcal{H}^{\aff})$-module of rank $|W|^2$.
\end{enumerate}
\end{Lemma}
The center of $\mathcal{H}^{\aff}$ also has a geometric interpretation:
\begin{equation*}
Z(\mathcal{H}^{\aff})   = \overline{\mathbb{Q}}_{\ell}[q,q^{-1}][X^*]^W \cong \mathcal{O}(T/W \times \mathbb{G}_m).
\end{equation*}
Hence, the central characters $\chi : Z(\mathcal{H}^{\aff}) \rightarrow \overline{\mathbb{Q}}_{\ell}$ are parametrized by the points of $T/W \times \mathbb{G}_m $ or equivalently by semisimple conjugacy classes in $G\times \mathbb{G}_m$. For $(s,q) \in G \times \mathbb{G}_m$ semisimple we denote the corresponding central character by $\chi_{(s,q)}$. By (the countable dimension version of) Schur's lemma \cite[Lemma 2.1.3]{chriss2009representation}, any simple $\mathcal{H}^{\aff}$-module admits a central character. Together with \cref{Lemma: algebraic properties of affine Hecke algebra}(4) this implies that the simple $\mathcal{H}^{\aff}$-modules are finite-dimensional. Moreover, the simple $\mathcal{H}^{\aff}$-modules with central character $\chi_{(s,q)}$ are precisely the simple modules of the finite-dimensional algebra
\begin{equation*}
\mathcal{H}^{\aff}_{(s,q)} :=  \mathcal{H}^{\aff} / (\ker(\chi_{(s,q)})).
\end{equation*}
Here $(\ker(\chi_{(s,q)})) = \mathcal{H}^{\aff} \cdot \ker(\chi_{(s,q)}) $ is the two-sided ideal in $\mathcal{H}^{\aff}$ generated by $ \ker(\chi_{(s,q)})$. The affine Hecke algebra also has a geometric incarnation which we recall next. Let $\mathfrak{g}$ (resp. $\mathfrak{b}$) be the Lie algebra of $G$ (resp. $B$), $\mathfrak{n} = [\mathfrak{b}, \mathfrak{b}]$ and $\mathcal{N} \subset \mathfrak{g}$ the nilpotent cone. The Springer resolution is the space
\begin{equation*}
\tilde{\mathcal{N}}:= G \times^B \mathfrak{n} \cong \{ (B', n) \in \mathcal{B} \times \mathcal{N} \mid n \in \mathfrak{b}' \}.
\end{equation*}
where we identify $\mathcal{B}$ with the set of all Borel subgroups $B'\subset G$ and denote by $\mathfrak{b}'$ the Lie algebra of $B'$. This comes with the two projections
\begin{equation}\label{eq: diagram for standard Springer resolution}
\begin{tikzcd}
& \tilde{\mathcal{N}} \arrow[ld, "\pi"'] \arrow[rd, "\mu"] & \\
\mathcal{B} & & \mathcal{N} \subset \mathfrak{g}.
\end{tikzcd}
\end{equation}
The maps $\pi$ and $\mu$ are $G\times\mathbb{G}_m$-equivariant where $t \in \mathbb{G}_m$ acts by scaling with $t^{-1}$ on $\mathcal{N}$ (resp. $\mathfrak{n}$) and trivially on $\mathcal{B}$. For any $G\times\mathbb{G}_m$-variety $X$, we denote by $ K^{G\times\mathbb{G}_m}(X)$ the equivariant algebraic $K$-theory of $X$ and we define
\begin{equation*}
K^{G\times\mathbb{G}_m}(X)_{\overline{\mathbb{Q}}_{\ell}} :=K^{G\times\mathbb{G}_m}(X)  \otimes_{\mathbb{Z}} \overline{\mathbb{Q}}_{\ell} .
\end{equation*}
Let $Z = \tilde{\mathcal{N}} \times_{\mathcal{N}} \tilde{\mathcal{N}}$ be the Steinberg variety. By \cite[(5.2.21)]{chriss2009representation} there is a convolution operation
\begin{equation*}
\star: K^{G \times \mathbb{G}_m} (Z) \otimes K^{G \times \mathbb{G}_m} (Z) \rightarrow K^{G \times \mathbb{G}_m} (Z)
\end{equation*}
which equips $K^{G \times \mathbb{G}_m} (Z)_{\overline{\mathbb{Q}}_{\ell}}$ with the structure of a $\overline{\mathbb{Q}}_{\ell}[q,q^{-1}]$-algebra (here multiplication by $q$ corresponds to tensoring with the irreducible weight $1$ representation of $\mathbb{G}_m$).
\begin{Theorem}\label{Theorem: Kazhdan-Lusztig isomorphism}\cite{kazhdan1987proof,chriss2009representation}
There is an isomorphism of $\overline{\mathbb{Q}}_{\ell}[q,q^{-1}]$-algebras
\begin{equation*}
\mathcal{H}^{\aff} \cong K^{G \times \mathbb{G}_m} ( Z)_{\overline{\mathbb{Q}}_{\ell}}.
\end{equation*}
Under this isomorphism, the center $Z(\mathcal{H}^{\aff})$ corresponds to $K^{G \times \mathbb{G}_m}(\pt)_{\overline{\mathbb{Q}}_{\ell}}$.
\end{Theorem}
The truncated Hecke algebra $\mathcal{H}^{\aff}_{(s,q)}$ also has a geometric description which we discuss next. Passing to $(s,q)$-fixed points in \eqref{eq: diagram for standard Springer resolution}, we obtain a diagram
\begin{equation*}
\begin{tikzcd}
& \tilde{\mathcal{N}}^{(s,q)} \arrow[ld, "\pi^{(s,q)}"'] \arrow[rd, "\mu^{(s,q)}"] & \\
\mathcal{B}^s & & \mathcal{N}^{(s,q)} \subset \mathfrak{g}^{(s,q)}.
\end{tikzcd}
\end{equation*}
Let
\begin{equation*}
\tilde{\mathcal{N}}^{(s,q)} = \bigsqcup_{i \in I} \tilde{\mathcal{N}}^{(s,q),i}
\end{equation*}
be the decomposition into connected components. The variety $\tilde{\mathcal{N}}^{(s,q)}$ is smooth (c.f. \cite[Lemma 5.11.1]{chriss2009representation}) and hence its connected components $\tilde{\mathcal{N}}^{(s,q),i}$ are also smooth. Thus, we can consider the constant perverse sheaves
\begin{equation*}
\mathcal{C}_{\tilde{\mathcal{N}}^{(s,q),i}} :=  \textbf{1}_{\tilde{\mathcal{N}}^{(s,q),i}} [\dim \tilde{\mathcal{N}}^{(s,q),i}] \in D^b_c(\tilde{\mathcal{N}}^{(s,q),i}, \overline{\mathbb{Q}}_{\ell})
\end{equation*}
and the corresponding $(s,q)$-Springer sheaves
\begin{align*}
\Sbf^{(s,q),i} &:= (\mu^{(s,q)})_* \mathcal{C}_{\tilde{\mathcal{N}}^{(s,q),i}} \in D^b_c(\mathcal{N}^{(s,q)}, \overline{\mathbb{Q}}_{\ell}) \\
\Sbf^{(s,q)} &:= \bigoplus_{i \in I} \Sbf^{(s,q),i}.
\end{align*}
\begin{Theorem}\label{Theorem: Ext version of Deligne-Langlands}\cite{chriss2009representation}
There is an isomorphism of $\overline{\mathbb{Q}}_{\ell}$-algebras
\begin{equation*}
\mathcal{H}^{\aff}_{(s,q)} \cong \Hom_{D^b_c(\mathcal{N}^{(s,q)}, \overline{\mathbb{Q}}_{\ell})}^*(\Sbf^{(s,q)}, \Sbf^{(s,q)}).
\end{equation*}
\end{Theorem}
\begin{proof}
This is proved in \cite[Proposition 8.1.5, Lemma 8.6.1]{chriss2009representation}. The only difference to our situation is that we work with the constructible derived category of $\overline{\mathbb{Q}}_{\ell}$-sheaves $D^b_c(X, \overline{\mathbb{Q}}_{\ell})$ coming from the (pro-)étale topology whereas \textit{loc.cit.} works with the constructible derived category $D^b_c(X(\mathbb{C}), \mathbb{C})$ coming from the associated complex analytic space $X(\mathbb{C})$. It turns out that the two approaches are equivalent: By \cite[p.146]{beilinson2018faisceaux} there is a fully faithful functor
\begin{equation*}
D^b_c(X, \overline{\mathbb{Q}}_{\ell}) \rightarrow D^b_c(X(\mathbb{C}), \mathbb{C})
\end{equation*}
(which involves choosing an isomorphism $\mathbb{C} \cong \overline{\mathbb{Q}}_{\ell}$). For $X = \mathcal{N}^{(s,q)}$ this identifies the $\overline{\mathbb{Q}}_{\ell}$-Springer sheaf $\Sbf^{(s,q)}$ with its analytic version and thus we get an isomorphism of graded algebras
\begin{equation*}
\Hom^*_{D^b_c(X, \overline{\mathbb{Q}}_{\ell})}(\Sbf^{(s,q)},\Sbf^{(s,q)}) \cong \Hom^*_{D^b_c(X(\mathbb{C}), \mathbb{C})}(\Sbf^{(s,q)},\Sbf^{(s,q)}).
\end{equation*}
However, since we have avoided analytic arguments so far, it is probably more naturally to prove the theorem directly in the $\overline{\mathbb{Q}}_{\ell}$-setting. It turns out that this can be done by essentially the same argument as in the analytic setting in \cite{chriss2009representation}: Let $A \subset G \times \mathbb{G}_m$ be the closed subgroup generated by $(s,q)$. Then $Z^{(s,q)} = Z^A$. Denote by $L_{(s,q)}$ the $1$-dimensional $Z(\mathcal{H}^{\aff})$-representation corresponding to $\chi_{(s,q)}$. Then there is a chain of algebra isomorphisms
\begin{align*}
\mathcal{H}^{\aff}_{(s,q)} & \cong K^{G \times \mathbb{G}_m}(Z)_{\overline{\mathbb{Q}}_{\ell}} \otimes_{K^{G \times \mathbb{G}_m}(\pt)_{\overline{\mathbb{Q}}_{\ell}}} L_{(s,q)}  \\
& \cong  K^{A}(Z)_{\overline{\mathbb{Q}}_{\ell}} \otimes_{K^{A}(\pt)_{\overline{\mathbb{Q}}_{\ell}}}  L_{(s,q)} \\
& \cong  K^{A}(Z^A)_{\overline{\mathbb{Q}}_{\ell}} \otimes_{K^{A}(\pt)_{\overline{\mathbb{Q}}_{\ell}}}  L_{(s,q)} \\
& \cong  K(Z^A)_{\overline{\mathbb{Q}}_{\ell}} \\
& \cong H_*(Z^{A}, \overline{\mathbb{Q}}_{\ell} ) \\
& =  H_*(Z^{(s,q)}, \overline{\mathbb{Q}}_{\ell} ) \\
& \cong \Hom_{D^b_c(\mathcal{N}^{(s,q)}, \overline{\mathbb{Q}}_{\ell})}^*(\Sbf^{(s,q)}, \Sbf^{(s,q)}).
\end{align*}
The first four algebra isomorphisms are exactly as in \cite[(8.1.6)]{chriss2009representation} (note that these are statements about equivariant algebraic $K$-theory of varieties over $\mathbb{C} \cong \overline{\mathbb{Q}}_{\ell}$, which does not involve the analytic topology). The last algebra isomorphism is proved exactly as \cite[Theorem 8.6.7]{chriss2009representation} which only relies on the six functor formalism. The remaining (fifth) isomorphism is given by the composition of the `Riemann-Roch map' for singular varieties and the cycle class map from \eqref{eq: cycle class map}:
\begin{equation}\label{eq: iso between K theory and homology}
K(Z^A)_{\overline{\mathbb{Q}}_{\ell}}  \overset{\sim}{\rightarrow}  A_*(Z^A)_{\overline{\mathbb{Q}}_{\ell}} \overset{\sim}{\rightarrow} H_*(Z^{A}, \overline{\mathbb{Q}}_{\ell} ).
\end{equation}
These maps are isomorphisms by \cite[Corollary 18.3.2]{fulton2013intersection} and \cref{Prop: H_*(Z) is spanned by fundamental classes} (together with \cref{Corollary: the sqSpringer resolution is of Springer type}). To prove that the isomorphism from \eqref{eq: iso between K theory and homology} is compatible with convolution, it suffices to show that the push, pull and $\otimes$ (resp. $\cap$) constructions that go into the definition of convolution are preserved. This can be found in \cite[Theorem 18.3]{fulton2013intersection} and \cite[Théorème 6.1,7.2]{laumon1976homologie}.
\end{proof}
As the notation suggests, the sheaf $\Sbf^{(s,q)}$ is a Springer sheaf in the sense of \cref{Section: morphisms of Springer type}. To check this, we need the following standard results about centralizers (see \cite[Proposition 8.8.7]{chriss2009representation}).
\begin{Lemma}\label{Lemma: fixed points of flag variety with respect to semisimple element}
\begin{enumerate}[label=(\roman*)]
\item The centralizer $G(s)$ is connected and reductive;
\item Each connected component of the fixed-point variety $\mathcal{B}^s$ is $G(s)$-equivariantly isomorphic to the flag variety of $G(s)$.
\end{enumerate}
\end{Lemma}
We pick the Borel subgroup $B \subset G$ such that $s \in B$, i.e. $B \in \mathcal{B}^s$. Note that $B(s) = G(s) \cap B$ is the stabilizer of $B \in \mathcal{B}^s$ in $G(s)$. By \cref{Lemma: fixed points of flag variety with respect to semisimple element}(ii) this implies that $B(s) \subset G(s)$ is a Borel subgroup. Let 
\begin{equation*}
\mathcal{B}^s = \bigsqcup_{j \in J} C_j
\end{equation*}
be the decomposition of $\mathcal{B}^s$ into connected components. Then for each $j \in J$ there is a unique element $B_j \in C_j \cong G(s)/B(s)$ whose stabilizer in $G(s)$ is $B(s)$ (i.e. $B_j(s) = B(s)$). Let $\mathfrak{b}_j$ be the Lie algebra of $B_j$ and $\mathfrak{n}_j = [\mathfrak{b}_j, \mathfrak{b}_j]$. The fiber of $B_j$ under $\pi^{(s,q)}$ is given by
\begin{equation*}
 (\pi^{(s,q)})^{-1} (B_j) \cong \mathcal{N}^{(s,q)} \cap \pi^{-1}(B_j) = \mathcal{N}^{(s,q)} \cap \mathfrak{b}_j = \mathfrak{n}_j^{(s,q)}.
\end{equation*}
Hence, we get an isomorphism (c.f. \cite[Corollary 8.8.9]{chriss2009representation})
\begin{equation*}
(\pi^{(s,q)})^{-1}( C_j ) \cong  G(s) \times^{B(s)} \mathfrak{n}_j^{(s,q)}.
\end{equation*}
Note that the variety $G(s) \times^{B(s)} \mathfrak{n}_j^{(s,q)}$ is connected, so the $(\pi^{(s,q)})^{-1}( C_j )$ are already the connected components $\tilde{\mathcal{N}}^{(s,q), i}$ of $\mathcal{N}^{(s,q)}$. We have thus shown that the $(s,q)$-Springer resolution is of Springer type (c.f. \cref{Definition: morphism of Springer type}).
\begin{Corollary}\label{Corollary: the sqSpringer resolution is of Springer type}
For each $i \in I$, the morphism $\mu^{(s,q)}: \tilde{\mathcal{N}}^{(s,q),i} \rightarrow \mathfrak{g}^{(s,q)}$ is of Springer type with Springer sheaf $\Sbf^{(s,q),i}$.
\end{Corollary}
Hence, we can consider the corresponding Springer category
\begin{equation*}
D_{\Spr}(\mathcal{N}^{(s,q)},\overline{\mathbb{Q}}_{\ell}) := D_{\Spr}(\mathfrak{g}^{(s,q)},\overline{\mathbb{Q}}_{\ell}) = \langle X_1, ..., X_n \rangle_{D^b_c(\mathfrak{g}^{(s,q)}, \overline{\mathbb{Q}}_{\ell}), \Delta}
\end{equation*}
where the $X_i$ are the simple perverse constituents of $\Sbf^{(s,q)}$.
\begin{Remark}
It might be more natural to consider the Springer category $D_{\Spr}(\mathcal{N}^{(s,q)},\overline{\mathbb{Q}}_{\ell}) $ as a subcategory of $D^b_c(\mathcal{N}^{(s,q)}, \overline{\mathbb{Q}}_{\ell})$ instead of $D^b_c(\mathfrak{g}^{(s,q)}, \overline{\mathbb{Q}}_{\ell})$. Note that this is not much of a difference since the canonical functor $D^b_c(\mathcal{N}^{(s,q)}, \overline{\mathbb{Q}}_{\ell}) \rightarrow D^b_c(\mathfrak{g}^{(s,q)}, \overline{\mathbb{Q}}_{\ell})$ induced by the closed immersion $\mathcal{N}^{(s,q)} \hookrightarrow \mathfrak{g}^{(s,q)}$ is fully faithful. In particular, we also have $D_{\Spr}(\mathcal{N}^{(s,q)},\overline{\mathbb{Q}}_{\ell}) = \langle X_1, ..., X_n \rangle_{D^b_c(\mathcal{N}^{(s,q)}, \overline{\mathbb{Q}}_{\ell}), \Delta}$.
\end{Remark}
By \cref{Corollary: Springer formality over C} and \cref{Corollary: the sqSpringer resolution is of Springer type} there is an equivalence  of triangulated categories
\begin{equation*}
D_{\Spr}(\mathcal{N}^{(s,q)}, \overline{\mathbb{Q}}_{\ell})^{\op} \cong D_{\perf}(\Hom^*(\textbf{S}^{(s,q)},\textbf{S}^{(s,q)}) - \dgMod)
\end{equation*}
which sends $\textbf{S}^{(s,q)}$ to the free dg-module $\Hom^*(\textbf{S}^{(s,q)},\textbf{S}^{(s,q)})$. Combining this with  \cref{Theorem: Ext version of Deligne-Langlands}, we arrive at the following `derived Deligne-Langlands correspondence'.
\begin{Theorem}\label{Theorem: derived Deligne-Langlands correspondence}
There is an equivalence of triangulated categories
\begin{equation*}
D_{\Spr}(\mathcal{N}^{(s,q)}, \overline{\mathbb{Q}}_{\ell})^{\op} \cong D_{\perf}(\mathcal{H}^{\aff}_{(s,q)}- \dgMod)
\end{equation*}
which sends $\Sbf^{(s,q)}$ to the free dg-module $\mathcal{H}^{\aff}_{(s,q)}$. Here we consider $\mathcal{H}^{\aff}_{(s,q)}$ as a dg-algebra with vanishing differential and grading induced by the $\Hom^*$-grading in \cref{Theorem: Ext version of Deligne-Langlands}.
\end{Theorem}

\bibliographystyle{alpha}
\bibliography{bibliography}

\begin{thebibliography}{BZCHN23}

\bibitem[AGV72]{SGA4}
Michael Artin, Alexander Grothendieck, and Jean-Louis Verdier.
\newblock {\em Theorie de Topos et Cohomologie Étale des Schemas. Tome II.
  (SGA4)}, volume 270 of {\em Lecture Notes in Mathematics}.
\newblock Springer, Berlin, 1972.

\bibitem[AMS18]{aubert2018graded}
Anne-Marie Aubert, Ahmed Moussaoui, and Maarten Solleveld.
\newblock Graded {H}ecke algebras for disconnected reductive groups.
\newblock In {\em Geometric aspects of the trace formula}, Simons Symposia,
  pages 23--84. Springer, 2018.

\bibitem[AR16]{achar2014modular}
Pramod Achar and Simon Riche.
\newblock Modular perverse sheaves on flag varieties. {II}: {Koszul} duality
  and formality.
\newblock {\em Duke Math. J.}, 165(1):161--215, 2016.

\bibitem[BBD82]{beilinson2018faisceaux}
Alexander Beilinson, Joseph Bernstein, and Pierre Deligne.
\newblock Faisceaux pervers.
\newblock {\em Ast{\'e}risque}, 100:5--171, 1982.

\bibitem[BD84]{bernstein1984centre}
Joseph Bernstein and Pierre Deligne.
\newblock Le “centre” de {B}ernstein.
\newblock In {\em Représentations des groupes reductifs sur un corps local},
  pages 1--32. Hermann, Paris, 1984.

\bibitem[Bei87]{beilinson1987derived}
Alexander Beilinson.
\newblock On the derived category of perverse sheaves.
\newblock In {\em K-theory, Arithmetic and Geometry}, volume 1289 of {\em
  Lecture Notes in Mathematics}, pages 27--41. Springer, Berlin, 1987.

\bibitem[BGS96]{beilinson1996koszul}
Alexander Beilinson, Victor Ginzburg, and Wolfgang Soergel.
\newblock Koszul duality patterns in representation theory.
\newblock {\em Journal of the American Mathematical Society}, 9(2):473--527,
  1996.

\bibitem[BK98]{bushnell1998smooth}
Colin Bushnell and Philip Kutzko.
\newblock Smooth representations of reductive $p$-adic groups: structure theory
  via types.
\newblock {\em Proceedings of the London Mathematical Society}, 77(3):582--634,
  1998.

\bibitem[BL94]{bernstein2006equivariant}
Joseph Bernstein and Valery Lunts.
\newblock {\em Equivariant sheaves and functors}, volume 1578 of {\em Lecture
  Notes in Mathematics}.
\newblock Springer, Berlin Heidelberg, 1994.

\bibitem[Bor76]{borel1976admissible}
Armand Borel.
\newblock Admissible representations of a semi-simple group over a local field
  with vectors fixed under an {I}wahori subgroup.
\newblock {\em Inventiones mathematicae}, 35(1):233--259, 1976.

\bibitem[BS01]{balmer2001idempotent}
Paul Balmer and Marco Schlichting.
\newblock Idempotent completion of triangulated categories.
\newblock {\em Journal of Algebra}, 236(2):819--834, 2001.

\bibitem[BS15]{bhatt2015pro}
Bhargav Bhatt and Peter Scholze.
\newblock The pro-étale topology for schemes.
\newblock {\em Ast{\'e}risque}, 369:99--201, 2015.

\bibitem[BZCHN20]{ben2020coherent}
David Ben-Zvi, Harrison Chen, David Helm, and David Nadler.
\newblock Coherent {S}pringer theory and the categorical {D}eligne-{L}anglands
  correspondence.
\newblock {\em \href{https://arxiv.org/abs/2010.02321}{arXiv:2010.02321}},
  2020.

\bibitem[BZCHN23]{ben2023between}
David Ben-Zvi, Harrison Chen, David Helm, and David Nadler.
\newblock Between coherent and constructible local {L}anglands correspondences.
\newblock {\em \href{https://arxiv.org/abs/2302.00039}{arXiv:2302.00039}},
  2023.

\bibitem[Cas80]{casselman1980unramified}
William Casselman.
\newblock The unramified principal series of $ p $-adic groups. {I}. {T}he
  spherical function.
\newblock {\em Compositio Mathematica}, 40(3):387--406, 1980.

\bibitem[CG97]{chriss2009representation}
Neil Chriss and Victor Ginzburg.
\newblock {\em Representation theory and complex geometry}.
\newblock Birkhäuser, Boston, 1997.

\bibitem[Cho20]{chough2020pro}
Chang-Yeon Chough.
\newblock The pro-{\'e}tale topology for algebraic stacks.
\newblock {\em Communications in Algebra}, 48(9):3761--3770, 2020.

\bibitem[DCLP88]{de1988homology}
Corrado De~Concini, George Lusztig, and Claudio Procesi.
\newblock Homology of the zero-set of a nilpotent vector field on a flag
  manifold.
\newblock {\em Journal of the American Mathematical Society}, 1(1):15--34,
  1988.

\bibitem[Del77]{SGA4-1/2}
Pierre Deligne.
\newblock {\em Cohomologie étale : Sèminaire de géométrie algébrique de
  Bois-Marie (SGA 4 1/2)}, volume 569 of {\em Lecture Notes in Mathematics}.
\newblock Springer, Berlin, 1977.

\bibitem[Del80]{deligne1980conjecture}
Pierre Deligne.
\newblock La conjecture de {W}eil: {II}.
\newblock {\em Publications Math{\'e}matiques de l'IH{\'E}S}, 52:137--252,
  1980.

\bibitem[DGMS75]{deligne1975real}
Pierre Deligne, Phillip Griffiths, John Morgan, and Dennis Sullivan.
\newblock Real homotopy theory of {K}{\"a}hler manifolds.
\newblock {\em Inventiones mathematicae}, 29:245--274, 1975.

\bibitem[Eke90]{ekedahl2007adic}
Torsten Ekedahl.
\newblock On the adic formalism.
\newblock In {\em The Grothendieck Festschrift}, pages 197--218. Birkhäuser,
  Boston, 1990.

\bibitem[ES22]{eberhardt2022motivic}
Jens~Niklas Eberhardt and Catharina Stroppel.
\newblock Motivic {S}pringer theory.
\newblock {\em Indagationes Mathematicae}, 33(1):190--217, 2022.

\bibitem[Ful98]{fulton2013intersection}
William Fulton.
\newblock {\em Intersection theory}.
\newblock Springer, New York, second edition, 1998.

\bibitem[Hel20]{hellmann2020derived}
Eugen Hellmann.
\newblock On the derived category of the {I}wahori-{H}ecke algebra.
\newblock {\em \href{https://arxiv.org/abs/2006.03013}{arXiv:2006.03013}},
  2020.

\bibitem[Jan87]{jantzen2003representations}
Jens Jantzen.
\newblock {\em Representations of algebraic groups}, volume 111 of {\em Pure
  and Applied Mathematics}.
\newblock Academic Press, Boston, 1987.

\bibitem[Kar78]{karoubi1978k}
Max Karoubi.
\newblock {\em K-theory: {A}n introduction}, volume 226 of {\em Grundlehren der
  mathematischen Wissenschaften}.
\newblock Springer, Berlin, 1978.

\bibitem[Kel94]{keller1994deriving}
Bernhard Keller.
\newblock Deriving {DG} categories.
\newblock {\em Annales scientifiques de l'Ecole normale sup{\'e}rieure},
  27(1):63--102, 1994.

\bibitem[KL87]{kazhdan1987proof}
David Kazhdan and George Lusztig.
\newblock Proof of the {D}eligne-{L}anglands conjecture for {H}ecke algebras.
\newblock {\em Inventiones mathematicae}, 87(1):153--215, 1987.

\bibitem[Lau76]{laumon1976homologie}
G{\'e}rard Laumon.
\newblock Homologie {\'e}tale.
\newblock {\em Astérisque}, 36-37:163--188, 1976.

\bibitem[LO08]{laszlo2008six}
Yves Laszlo and Martin Olsson.
\newblock The six operations for sheaves on {A}rtin stacks {II}: adic
  coefficients.
\newblock {\em Publications Math{\'e}matiques de l'IH{\'E}S}, 107:169--210,
  2008.

\bibitem[Lus89]{lusztig1989affine}
George Lusztig.
\newblock Affine {H}ecke algebras and their graded version.
\newblock {\em Journal of the American Mathematical Society}, 2(3):599--635,
  1989.

\bibitem[Lus95a]{lusztig1995classification}
George Lusztig.
\newblock Classification of unipotent representations of simple $p$-adic
  groups.
\newblock {\em International Mathematics Research Notices}, 1995(11):517--589,
  1995.

\bibitem[Lus95b]{lusztig1995cuspidal}
George Lusztig.
\newblock Cuspidal local systems and graded {H}ecke algebras. {II}.
\newblock In {\em Representations of groups}, volume~16 of {\em Conference
  Proceedings, Canadian Mathematical Society}, pages 217--275, 1995.

\bibitem[Nee01]{neeman2001triangulated}
Amnon Neeman.
\newblock {\em Triangulated categories}, volume 148 of {\em Annals of
  mathematics studies}.
\newblock Princeton University Press, Princeton, 2001.

\bibitem[PVdB19]{polishchuk2019semiorthogonal}
Alexander Polishchuk and Michel Van~den Bergh.
\newblock Semiorthogonal decompositions of the categories of equivariant
  coherent sheaves for some reflection groups.
\newblock {\em Journal of the European Mathematical Society}, 21(9):2653--2749,
  2019.

\bibitem[Rid13]{rider2013formality}
Laura Rider.
\newblock Formality for the nilpotent cone and a derived {S}pringer
  correspondence.
\newblock {\em Advances in Mathematics}, 235:208--236, 2013.

\bibitem[RR21]{rider2021formality}
Laura Rider and Amber Russell.
\newblock Formality and {L}usztig’s generalized {S}pringer correspondence.
\newblock {\em Algebras and Representation Theory}, 24:699--714, 2021.

\bibitem[RSW14]{riche2014modular}
Simon Riche, Wolfgang Soergel, and Geordie Williamson.
\newblock Modular {K}oszul duality.
\newblock {\em Compositio Mathematica}, 150(2):273--332, 2014.

\bibitem[Sch11]{schnurer2011equivariant}
Olaf Schn{\"u}rer.
\newblock Equivariant sheaves on flag varieties.
\newblock {\em Mathematische Zeitschrift}, 267(1-2):27--80, 2011.

\bibitem[Sol22]{solleveld2022graded}
Maarten Solleveld.
\newblock Graded {H}ecke algebras and equivariant constructible sheaves on the
  nilpotent cone.
\newblock {\em \href{https://arxiv.org/abs/2205.07490}{arXiv:2205.07490}},
  2022.

\bibitem[{Sta}23]{stacks-project}
The {Stacks Project Authors}.
\newblock \textit{Stacks Project}.
\newblock \url{https://stacks.math.columbia.edu}, 2023.

\end{thebibliography}

\end{document}